\documentclass[12pt]{elsarticle}
\usepackage[top=2cm,bottom=2cm,left=2.5cm,right=2.5cm,marginparwidth=1.75cm]{geometry}
\usepackage[english]{babel}

\usepackage{textcomp}
\usepackage[super]{nth}
\usepackage{float}
\usepackage{subfloat}
\usepackage{subcaption}
\usepackage{multirow}
\usepackage{makecell}
\usepackage{longtable,tabularx}
\usepackage{lipsum}
\usepackage{pgfplots}
\usepackage{tikz}
\usetikzlibrary{positioning,external}
\usepackage{algorithm}
\usepackage{algpseudocode}
\usepackage{amsmath}
\usepackage{pifont}% http://ctan.org/pkg/pifont% \usepackage{amssymb}
\usepackage{amstext}
\usepackage{amsfonts}
\usepackage{bm}
\usepackage{mathtools}
\usepackage{amssymb}
\usepackage{comment}
\usepackage{acronym} % For acronyms
\usepackage{hyperref}
\hypersetup{
    colorlinks=false,
    linkcolor=black,
    filecolor=magenta,      
    urlcolor=cyan
    }

\usepackage[bb=boondox]{mathalfa}

\usepackage[capitalize]{cleveref}
\usepackage{apptools}
\usepackage{graphicx,import}
\usepackage{lastpage}
\usepackage{fancyhdr}
\usepackage{ifthen}
\usepackage[]{siunitx}
\usepackage{algorithm}
\usepackage{algpseudocode}

\newcommand*{\doublerule}{\hrule width \hsize height 1pt \kern 0.5mm \hrule width \hsize height 2pt}
\newcommand\doublerulefill{\leavevmode\leaders\vbox{\hrule width .1pt\kern1pt\hrule}\hfill\kern0pt }
\usetikzlibrary{plotmarks}
\pgfplotsset{compat=newest,
every plot/.append style={color=black, mark=none},
every axis/.append style={
label style={font=\fontsize{9pt}{1em}},
tick style={font=\fontsize{9pt}{1em}\selectfont, color=black, line cap=round},
ticklabel style= {font=\fontsize{9pt}{1em}\selectfont},
legend style={legend cell align=left, font=\fontsize{9pt}{1em}\selectfont, align=left},
/pgf/number format/1000 sep={},
yticklabel style={
        /pgf/number format/fixed,
        /pgf/number format/precision=5
},
scaled y ticks=false
}}
\newcommand{\ipc}{P_{IC}}

\newcommand{\pc}{P_C}

\newcommand{\pclim}{\bar{P}_C}
\newcommand{\dmiss}{d_{miss}}

\newcommand{\dv}{\Delta v}
\newcommand{\cmark}{\ding{51}}%
\newcommand{\xmark}{\ding{55}}%

\let\oldhat\hat
\renewcommand{\vec}[1]{\boldsymbol{#1}}
\renewcommand{\hat}[1]{\oldhat{\boldsymbol{#1}}}

\acrodef{2PBVP}{Two-Point Boundary Value Problem} 
\acrodef{3PBVP}{Three-Point Boundary Value Problem} 
\acrodef{ACPL}{Accepted Collision Probability Level} 
\acrodef{AI}{Artificial Intelligence}
\acrodef{AIDA}{Accurate Integrator for Debris Analysis}
\acrodef{ADR}{Active Debris Removal}
\acrodef{AOCS}{Attitude and Orbit Control System} 
\acrodef{ARES}{Assessment of Risk Event Statistics} 
\acrodef{AvA}{Active vs Active}
\acrodef{BVP}{Boundary Value Problem}
\acrodef{BP}{B-plane}
\acrodef{CAOIA}{Conjunction Assessment Operations Implementation Agreement} 
\acrodef{CA}{Conjunction Analysis}
\acrodef{CAM}{Collision Avoidance Manoeuvre} 
\acrodef{CARA}{Conjunction Risk Assessment}
\acrodef{CONAN}{Conjunction Analysis}
\acrodef{COLA}{Collision Avoidance}
\acrodef{CDM}{Conjunction Data Message}
\acrodef{CEGAR}{Counter-Example Guided Abstraction Refinement} 
\acrodef{HFEM}{higher-Fidelity Ephemeris Model}
\acrodef{CPU}{Central Processing Unit} 
\acrodef{CR3BP}{Circular Restricted Three Body Problem}
\acrodef{CW}{Clohessy-Wiltshire}
\acrodef{DA}{Differential Algebra}
\acrodef{IAC}{International Astronautical Congress}
\acrodef{DELTA}{Debris Environment Long-Term Analysis} 
\acrodef{DISCOS}{Database and Information System Characterizing Objects in Space} 
\acrodef{DRAMA}{Debris Risk Assessment and Mitigation Analysis}
\acrodef{ECI}{Earth-Centered Inertial} 
\acrodef{ECEF}{Earth-Centered Earth-Fixed}
\acrodef{EO}{Energy-Optimal} 
\acrodef{EOP}{Energy-Optimal Problem} 
\acrodef{EOCP}{Energy-Optimal Control Problem} 
\acrodef{AIAA}{American Institute of Aeronautics and Astronautics} 
\acrodef{EOE}{Equinoctial Orbital Elements} 
\acrodef{ESA}{European Space Agency} 
\acrodef{FBC}{Finite Burn Conversion} 
\acrodef{FO}{Fuel-Optimal} 
\acrodef{FOP}{Fuel-Optimal Problem} 
\acrodef{FPE}{Fokker-Planck Equation}
\acrodef{GEO}{Geostationary Orbit}
\acrodef{GMM}{Gaussian Mixture Model}
\acrodef{GMMs}{Gaussian Mixture Models}
\acrodef{GA}{Genetic Algorithm}
\acrodef{GNC}{Guidance, Navigation, and Control}
\acrodef{GSOC}{German Space Operation Center}
\acrodef{GTO}{Geostationary Transfer Orbit}
\acrodef{HBR}{Hard Body Radius}
\acrodef{IADC}{Inter-Agency Space Debris Coordination Committee}
\acrodef{IAF}{International Astronautical Federation}
\acrodef{IRA}{Iterative Risk Allocation}
\acrodef{IPoC}{Instantaneous Probability of Collision}
\acrodef{COLA}{COLlision Avoidance}
\acrodef{ISS}{International Space Station}
\acrodef{JAXA}{Japan Aerospace eXploration Agency}
\acrodef{JSpOC}{Joint Space Operations Center}
\acrodef{KDE}{Kernel Density Estimator}
\acrodef{KIAM}{Keldysh Institute of Applied Mathematics}
\acrodef{KKT}{Karush-Kuhn-Tucker}
\acrodef{KOZ}{Keep-Out-Zone}
\acrodef{LEGEND}{LEO-to-GEO Environment Debris model}
\acrodef{LEO}{Low Earth Orbit} 
\acrodef{LEOMEGCONST}{Low Earth Orbit Mega-Constellation Transfers}
\acrodef{LLO}{Low Lunar Orbit}
\acrodef{LNT}{Lethal Non-Trackable} 
\acrodef{LVLH}{Local Vertical Local Horizon}
\acrodef{MC}{Monte Carlo}
\acrodef{MEE}{Modified Equinoctial Elements}
\acrodef{MEO}{Medium Earth Orbit}
\acrodef{MECSA}{Method of Equivalent Cross-Section}
\acrodef{MGMMs}{Multidirectional Gaussian Mixture Models}
\acrodef{MILP}{Mixed Integer Linear Program}
\acrodef{ML}{Machine Learning}
\acrodef{MM}{Monodromy Matrix}
\acrodef{MOCAT}{MIT Orbital Capacity Assessment Tool}
\acrodef{MOID}{Minimum Orbit Intersection Distance}
\acrodef{MASTER}{Meteoroid and Space Debris Terrestrial Environment Reference} 
\acrodef{MD}{Miss Distance}
\acrodef{NEODEEM}{Near-Earth Orbit Debris Environment Evolutionary Model}
\acrodef{NASA}{National Aeronautics and Space Administration}
\acrodef{NLI}{Non-Linearity Index}
\acrodef{NFR}{No-Further-Release}
\acrodef{NLP}{Non-Linear Program}
\acrodef{NRHO}{Near-Rectilinear Halo Orbits} 
\acrodef{OCP}{Optimal Control Problem}
\acrodef{OD}{Orbit Determination}
\acrodef{ODE}{Ordinary Differential Equation}
\acrodef{ODMSP}{Orbital Debris Mitigation Standard Practices}
\acrodef{ODPO}{Orbital Debris Program Office}
\acrodef{OPM}{Orbit Parameter Message}
\acrodef{ORDEM}{Orbital Debris Engineering Model}
\acrodef{ORSK}{Orbit Raising and Station Keeping}
\acrodef{PC}{Polynomial Chaos}
\acrodef{PDF}{Probability Density Function}
\acrodef{PMD}{Post Mission Disposal}
\acrodef{PoC}{Probability of Collision}
\acrodef{PTO}{Point-to-Orbit}
\acrodef{PTP}{Point-to-Point}
\acrodef{QCQP}{Quadratically Constrained Quadratic Program}
\acrodef{QP}{Quadratic Program}
\acrodef{RAE}{Royal Aircraft Establishment}
\acrodef{RAAN}{Right Ascension of the Ascending Node}
\acrodef{RL}{Reinforcement Learning}
\acrodef{RK}{Runge-Kutta}
\acrodef{ROD}{Refined Orbit Determination} 
\acrodef{r.f.}{reference frame} 
\acrodef{RSO}{Resident Space Object}
\acrodef{RTN}{Radial-Transverse-Normal}
\acrodef{RV}{Random Variable} 
\acrodef{S3TOC}{Spanish SST Operations Centre}
\acrodef{SCP}{Sequential Convex Program}
\acrodef{SCvx}{Successive Convexification}
\acrodef{SDA}{Space Data Association}
\acrodef{SDP}{Semidefinite Program}
\acrodef{SEM}{Space Environment Management}
\acrodef{SF}{Solar Flux}
\acrodef{SGP4}{Simplified General Perturbation 4}
\acrodef{SK}{Station Keeping}
\acrodef{SMD}{Squared Mahalanobis Distance}
\acrodef{SOCP}{Second-Order Cone Program}
\acrodef{SRP}{Solar Radiation Pressure}
\acrodef{SSA}{Space Situational Awareness}
\acrodef{SSO}{Sun-Synchronous Orbit}
\acrodef{SSN}{Space Surveillance Network}
\acrodef{SST}{Space Surveillance and Tracking}
\acrodef{STTs}{State Transition Tensors}
\acrodef{STM}{Space Traffic Management}
\acrodef{TCA}{Time of Closest Approach}
\acrodef{TFC}{Theory of Functional Connections} 
\acrodef{TIPoC}{Total Instantaneous Probability of Collision}
\acrodef{TLE}{Two-Line Element}
\acrodef{TRR}{Trust Region Radius}
\acrodef{UP}{Uncertainty Propagation} 
\acrodef{USSTRATCOM}{United States Strategic Command}
\acrodef{UT}{Unscented Transform} 
\acrodef{UTC}{Universal Time Coordinated}
\acrodef{ZOH}{Zeroth-Order Hold}
\acrodef{LP}{Linear Program}
\acrodef{TPoC}{Total Probability of Collision}
\acrodef{PP}{Polynomial Program}
\acrodef{VLEO}{Very Low Earth Orbit}
\acrodef{CORAM}{Collision Risk Assessment and Avoidance Manoeuvres}

\journal{Acta Astronautica}

\begin{document}
\sloppy
\begin{frontmatter}

\title{CAMmary: A Review of Spacecraft Collision Avoidance Manoeuvre Design Methods}

% Group authors per affiliation:
\author[Auckland]{Zeno Pavanello\corref{mycorrespondingauthor}}
\ead{zpav176@aucklanduni.ac.nz}
\author[Milano]{Luigi De Maria}
\author[Milano]{Andrea De Vittori}
\author[Milano]{Michele Maestrini}
\author[Milano]{Pierluigi Di Lizia}
\author[Auckland]{Roberto Armellin}

\address[Auckland]{Te P\=unaha \=Atea - The Space Institute, The University of Auckland, 20 Symonds Street, Auckland, New Zealand, 1010}
\address[Milano]{Department of Aerospace Science and Technology (DAER), Politecnico di Milano, Piazza Leonardo da Vinci 32, Milano, Itly, 20133}

\cortext[mycorrespondingauthor]{Corresponding author}

\begin{abstract}
    Ensuring safety for spacecraft operations has become a paramount concern due to the proliferation of space debris and the saturation of valuable orbital regimes. In this regard, the \ac{CAM} has emerged as a critical requirement for spacecraft operators, aiming to efficiently navigate through potentially hazardous encounters. Currently, when a conjunction is predicted, operators dedicate a considerable amount of time and resources to designing a \ac{CAM}. Given the increased frequency of conjunctions, autonomous computation of fuel-efficient \ac{CAM}s is crucial to reduce costs and improve the performance of future operations. 
To facilitate the transition to an autonomous \ac{CAM} design, it is useful to provide an overview of its state-of-the-art. In this survey article, a collection of the most relevant research contributions in the field is presented. We review and categorize existing \ac{CAM} techniques based on their underlying principles, such as (i) analytic, semi-analytic, or numerical solutions; (ii) impulsive or continuous thrust; (iii) deterministic or stochastic approaches, (iv) free or fixed manoeuvring time; (v) free or fixed thrust direction. Finally, to determine the validity of the algorithms potentially implementable for autonomous use, we perform a numerical comparison on a large set of conjunctions. With this analysis, the algorithms are evaluated in terms of computational efficiency, accuracy, and optimality of the computed policy. Through this comprehensive survey, we aim to provide insights into the state-of-the-art \ac{CAM} methodologies, identify gaps in current research, and outline potential directions for future developments in ensuring the safety and sustainability of spacecraft operations in increasingly congested orbital environments. 
\end{abstract}

\begin{keyword}
Survey\sep Collision Avoidance Manoeuvre\sep Trajectory Optimization\sep Space Traffic Management
\end{keyword}

\end{frontmatter}

\acresetall
\section{Introduction}
\label{sec:overallintroduction}
Since the dawn of the space age on October 4, 1957, space debris has consistently outnumbered operational satellites in Earth's orbit. Identified as a concern in the 1960s and gaining international attention by the 1970s, this growing collection of non-functional objects poses a significant challenge for spacefaring nations. In 1978, Donald Kessler \cite{Kessler_syn} postulated that collisions and explosions in orbit could lead to an exponential increase in debris, eventually rendering spaceflight too hazardous.

The 2024 \ac{ESA} Annual Space Environment Report \cite{ESA2024}, based on estimates from the MASTER model \cite{Horstmann2019}, provides the following space debris statistics: (i) Approximately 35,000 objects larger than 10~\si{cm} are typically monitored by ground-based sensors; (ii) around 900,000 objects ranging from 1~\si{cm} to 10~\si{cm} are challenging to detect yet sufficiently large to threaten operational satellites; (iii) an estimated 128 million objects between 1~\si{mm} and 1~\si{cm} are generally less hazardous as satellite shielding can mitigate their impact.
The \ac{LEO} and \ac{GEO} regions are particularly concerning, as they are becoming increasingly crowded due to the rising number of operational satellites, including micro-satellites and large constellations like Starlink.

Among the international community, growing concerns about the potentially catastrophic consequences of an unchecked evolution of the space population have led to the development of the broader concepts of \ac{SSA}\footnote{\url{https://www.euspa.europa.eu/eu-space-programme/ssa} Accessed on December 16 2024.} -- which encompasses \ac{SST} -- and \ac{STM}.\footnote{\url{https://defence-industry-space.ec.europa.eu/eu-space/space-traffic-management_en}, Accessed on November 25 2024}
The main elements that characterize \ac{STM} practices are space debris \textit{mitigation}
and \textit{remediation}. 
While remediation aims at actively reducing the number of space debris in orbit, mitigation involves measures to prevent the creation of new space debris. Common mitigation strategies include the intelligent design of satellite and rocket bodies, the implementation of end-of-life plans, and the actuation of \ac{CA} and \ac{COLA} routines when a close approach is foreseen.
In the aftermath of the 2009 collision between Iridium 33 and COSMOS 2251 \cite{KosmosIridium}, which highlighted deficiencies in collision assessment and communication, coordination among all entities involved in \ac{SSA} was initiated.
The UN \cite{UNGUIDELINES}, the \ac{IADC},\footnote{\url{https://iadc-home.org/what_iadc}, Accessed on 25 November 2024} the \ac{ODMSP},\footnote{\url{https://orbitaldebris.jsc.nasa.gov/library/usg_orbital_debris_mitigation_standard_practices_november_2019.pdf}. Accessed on December 19 2024.} and the \cite{ESASpaceDebrisMitigationWorkingGroup2023}
have been crucial promoters of international cooperation on space
debris mitigation, providing common guidelines to global space mission design standards. Moreover, since 2004, \ac{ESA} has provided \ac{DRAMA}, a specialized software tool for the compliance analysis of a space mission with space debris mitigation standards \cite{DRAMA}.

Reducing the risk of a collision fundamentally breaks down into two consecutive practices: \ac{CA} and \ac{COLA}, which, according to the \ac{AIAA}, are the top priority for the success of \ac{STM} operations \cite{AIAA2017}.
The conjunction screening process serves as a critical preliminary gatekeeper to \ac{CA}, identifying potentially dangerous conjunctions for further analysis. The earliest and most extensive screening process for all orbiting objects is carried out by the U.S. Space Force, which performs All-vs-All propagations every 24 hours to check for potential close encounters \cite{Krage2023}. When a close approach is detected, a \ac{CDM} \cite{CDM} is included in the Space-Track catalogue\footnote{\url{https://www.space-track.org}, Accessed on January 6 2025.} and the operator(s) of the involved spacecraft is/are notified. \ac{CDM}s serve as proximity alerts rather than collision warnings, requiring further risk assessment to determine if a collision is likely after \ac{CA} screening. The data in a \ac{CDM} allows spacecraft operators to make informed decisions about whether a \ac{CAM} is necessary: they typically review the information, perform further analysis – such as a refined orbit determination action on their satellite – and decide whether to manoeuvre.

Space agencies and operators have developed specific tools tailored to address the need for \ac{CAM} design. Since 2004, \ac{NASA} has used \ac {CARA} to process ephemeris data to identify close approaches and assess the risk of a collision \cite{CARA}; moreover, the software can help on-ground \ac{CAM} planning and
execution to mitigate dangerous events identified in the \ac{CA} phase. European counterparts to \ac{CARA} include the \ac{GSOC} pipeline \cite{aida2010collision}, the software Conjunction Analysis and
Evaluation Service: Alerts and Recommendations (CAESAR) from \textit{Centre National d’Études Spatiales} (CNES), the \ac{CONAN} tool from the EU \ac{SST} consortium \cite{CONAN}, and \ac{CORAM} from the \ac{ESA} Space Debris Office \cite{SanchezOrtiz2014}.
The actual process of analysing the \ac{CDM}, planning, uploading, and executing the manoeuvre has always been performed manually by operators \cite{Cortesi2024}. Even with the aid of the previously mentioned tools, this process requires a large number of work hours and 24/7 availability of specialized experts \cite{Flohrer2022}.

In this context, the rise of satellite mega-constellations like Starlink has increased the frequency of \ac{CAM}s. SpaceX, for instance, employs a conservative manoeuvre threshold of $10^{-5}$ on \ac{PoC}, averaging approximately 12 manoeuvres per satellite yearly \cite{Goldman2022}. Analysis of these constellations reveals a strong correlation between manoeuvre frequency and risk reduction.
A human-in-the-loop approach to the \ac{CA}/\ac{COLA} process will soon become unsustainable from an economical and practical standpoint. Introducing autonomous decision-making and actuation will become, therefore, fundamental in shaping future operational pipelines. This might be the main reason for the significant rise in works related to \ac{CAM} presented at the largest aerospace conference, the \ac{IAC}. The number of \ac{CAM}-focused contributions increased from 36 in 2021 to 56 in 2023 and further to 93 in 2024.

Autonomous operations are characterized by the ability to perform tasks and make decisions independently without real-time human intervention.
A key distinction can be drawn between \textit{ground-based autonomous} systems and \textit{onboard autonomous} ones. In the former, the whole \ac{CDM}-reading, \ac{CA}, and \ac{COLA} process is managed by a machine on the ground; to be fully autonomous, the system must also be able to upload the maneuvering command to the spacecraft at a suitable time, i.e., during a passage of the spacecraft above the ground station. Onboard autonomous systems, on the other hand, would realistically only require a single input from the ground, the \ac{CDM} upload, and they would manage the whole cascaded \ac{CA}/\ac{COLA} on board. Therefore, the key difference lies in the location of the decision-making process. Ground-based computers have exponentially higher computational power than onboard ones, so they can run more complex algorithms; however, being able to run the process on board would allow for timely adjustments of the maneuver in the case of unforeseen occurrences.

To this day, reliance on ground operations is essential for the extensive conjunction screening, the \ac{CDM} generation process, and the \ac{CA} process. Therefore, the only section of the whole \ac{CA}/\ac{COLA} pipeline that would realistically be implementable on the onboard computer of a spacecraft is the \ac{CAM} design/optimisation, devoid of any decision-making on whether or not to actuate the manoeuvre. However, since this last part of the whole pipeline is kick-started by an input coming from the ground station (\ac{CDM}/\ac{CA}), one could argue that it is more efficient to perform the optimisation of the manoeuvre itself on the ground, uploading it when possible. For this reason, fully ground-based systems are still the best option for an autonomous pipeline. If spacecraft had the ability to autonomously predict a close encounter and analyse the conjunction in real-time, the onboard option would be preferable, but this will realistically not be the case for the foreseeable future.

\ac{ML} is a generally accepted solution for the introduction of (on-ground) autonomous \ac{CA}/decision-making \cite{Gonzalo2020y,Gonzalo2021y,Sanchez2021,Stroe2021,Henry2023}, but classic optimisation techniques are still preferable for \ac{CAM} design.  
In this regard, the transition to autonomous \ac{COLA} requires mainly four characteristics of manoeuvre design: 
\begin{itemize}
    \item \textit{Optimality}: ensuring the computed manoeuvre is fuel-optimal is crucial. Fuel is a limited resource in space missions, and optimising its use can extend the mission's lifespan. Therefore, optimal control theory should be used to achieve at least locally optimal manoeuvres.
    \item \textit{Efficiency}: light and efficient algorithms reduce computational load and power consumption, which are critical factors in spacecraft operations.
    \item \textit{Reliability}: the algorithm must be guaranteed to produce the expected outcome under nominal conditions.
    \item \textit{Robustness}: spacecraft operate in environments where the modelling uncertainty is high, so the ability to handle off-nominal conditions, such as thrust errors and unexpected environmental perturbations, is pivotal to improve the chance of a successful manoeuvre.
\end{itemize}

The shift toward autonomous systems for \ac{STM} is of interest to many players in the field, which are developing projects to this end.
Key examples are Starling, involving SpaceX's Starlink constellation and \ac{NASA} \cite{Cannon2023Starling} and Collision Risk Estimation and Automated Mitigation (CREAM), from \ac{ESA}. The former aims to demonstrate advanced onboard \ac{CA} capabilities, such as continuous monitoring of both passive and active manoeuvring objects, along with planned manoeuvres, minimizing the need for direct operator-to-operator communication. The latter, framed in
the broader European \textit{Space Safety Programme} \cite{BastidaVirgili2019}, has the objective of demonstrating \lq\lq [...] a ‘decision support system’ for \ac{CAM}s, the integration and deployment of a pilot for automated and secure exchange of spacecraft manoeuvre plans and enhancement of autonomous decisions.”

Given the growing need for \ac{CAM} design, the goal of this work is to provide a comprehensive review of current state-of-the-art solutions, offering an in-depth analysis of existing models and a numerical comparison of those that are deemed applicable to autonomous use.
The paper is structured into seven sections. \cref{sec:PoC} introduces the tools necessary to better understand all steps required for \ac{CAM} design and their limitations, such as collision probability estimation and state uncertainties propagation. \cref{sec:NumericalCAM} discusses numerical approaches to \ac{CAM} design, including direct optimisation, heuristics, and \ac{AI}-based methods. Subsequently, \cref{sec:AnalyticalCAM} examines analytical and semi-analytical optimisation techniques. \cref{sec:Comparison} compares selected algorithms in terms of accuracy and computational efficiency. Finally, this paper explores the future trends in \ac{COLA} in \cref{sec:FutureTrends}, followed by the conclusions in \cref{sec:Conclusions}.

%%%%%%%%%%%%%%%%%%%%%%%%%%%%%%%%%%%%%%%%%%%%%%%%%%%%%%%%%%%%%%%%%%%%%%%%%%%
% PoC
%%%%%%%%%%%%%%%%%%%%%%%%%%%%%%%%%%%%%%%%%%%%%%%%%%%%%%%%%%%%%%%%%%%%%%%%%%%
% \section{Computational Techniques for Probability of Collision}
\section{Probability of Collision}
\label{sec:PoC}

This section discusses uncertainty propagation techniques and their characteristics and safety metrics for different types of close encounters.

%%%%%%%%%%%%%%%%%%%%%%%%%%%%%%%%%%%%%%%%%%%%%%%%%%%%%%%%%%%%%%%%%%
\subsection{Uncertainty Propagation}
\label{subsec:UP}

In \ac{CA} procedures, uncertainty propagation is essential for determining key satellite statistical information, such as the mean, covariance matrix, and \ac{PDF} of the state. Orbital dynamics with uncertainty can be modelled by the Itô stochastic differential equation \cite{maybeck1982_book}. The \ac{PDF}  time evolution for such a system is governed by the \ac{FPE} \cite{Fuller_art}. Solving it in orbital mechanics is challenging due to the high-dimensional state space and non-linear dynamics involved. To date, no satisfactory solution that can be deployed in real scenarios has been found \cite{Kumar} \cite{Sun_prop}.

Another rigorous approach to obtain a comprehensive statistical description involves employing \ac{MC} models. These simulations offer means to handle non-linear and non-Gaussian uncertainty propagation \cite{Junkins1996NonGaussianEP_article}. As \ac{PoC} decreases, accurately estimating it with \ac{MC} methods requires an increasing number of samples, significantly raising computational demands \cite{dagumMC}. Methods like Subset Simulation and Importance Sampling \cite{morselli2015high} reduce sample requirements but are less robust than standard \ac{MC} analyses.
In recent years, linear propagation techniques have been developed to balance accuracy and computational burden. In linear uncertainty propagation methods, Gaussian uncertainties are commonly assumed so that only the propagated mean value and the covariance matrix are needed to capture the entire \ac{PDF} at any time. In addition, local linearisation enables the linear propagation of the mean and covariance matrix using the state transition matrix, derived from the solution of the linearised system (LinCov) \cite{Battin1987AnIT,RAO}. On the other hand, statistical linearization approximates the dynamics in an `average' sense, allowing for the derivation of statistical differential equations for the mean and covariance matrix evolution \cite{gelb1974}. Despite their simplicity and efficiency, these methods become inaccurate with highly non-linear systems, long propagation windows, and high initial uncertainty.

Various methods for nonlinear uncertainty propagation have been developed as efficient alternatives to \ac{MC} simulations. \ac{UT} provides a second-order approximation of the first two moments of a Gaussian distribution by non-linearly integrating a few deterministically selected `sigma points' from the initial distribution, which may not suffice for long simulations \cite{Julier}. Similarly, cubature methods use deterministic sampling and integration, ensuring exact integral approximations for specific test functions, with spherical-radial cubature rules being particularly popular for their speed and robustness \cite{Arasaratnam}. \ac{PC} manages inputs and outputs via series approximations involving standard random variables, with some formulations allowing the use of existing orbit propagation tools as black-box components to provide accurate higher-order moment information. However, the number of \ac{PC} terms grows exponentially with input uncertainty dimensionality \cite{Arasaratnam}. Park and Scheeres \cite{Park} propose \ac{STTs} as a semi-analytic technique for orbit uncertainty propagation, but the variational approach used for their computation can suffer limitations for high-fidelity dynamic systems. \ac{DA} techniques have been proposed to overcome the limitation of the variational approach, enabling accurate non-linear uncertainty propagators with arbitrary order expansions \cite{Berz1999ModernMM_PDF}. \ac{DA} has successfully been used to reduce the computational effort of nonlinear uncertainty propagation \cite{Massari2017}, even in combination with \ac{GMMs} \cite{Sun2019} which approximate arbitrary \ac{PDF}s using a finite sum of weighted Gaussian kernels. Although flexible, they suffer from the curse of dimensionality and lack adaptive algorithms for weight selection. Various schemes, including \ac{KDE} \cite{Crisan} and adaptive \ac{GMMs} with dynamic splitting rules \cite{Vittaldev2016Space}, exist. Each strategy has its advantages and drawbacks, which makes them suitable for different scenarios, depending on the specific requirements and constraints of the problem \cite{maestrini2023electrocam}.

%%%%%%%%%%%%%%%%%%%%%%%%%%%%%%%%%%%%%%%%%%%%%%%%%%%%%%%%%%%%%%%%%%
\subsection{Probability of Collision Models}
\label{subsec:PoCModels}
Once the state uncertainties are characterised at each time of the conjunction assessment analyses, they can be used for \ac{PoC} estimation.
\ac{PoC} is the most comprehensive safety metric used in satellite operations. This parameter rigorously integrates \ac{MD}, covariance size and orientation, and the dimensions of the objects involved, making it comparable to probabilities of other failure scenarios, such as thruster failures \cite{Alfano2022_SatOperatorsSurvey}.

\begin{figure}[tb!]
    \centering
    \input{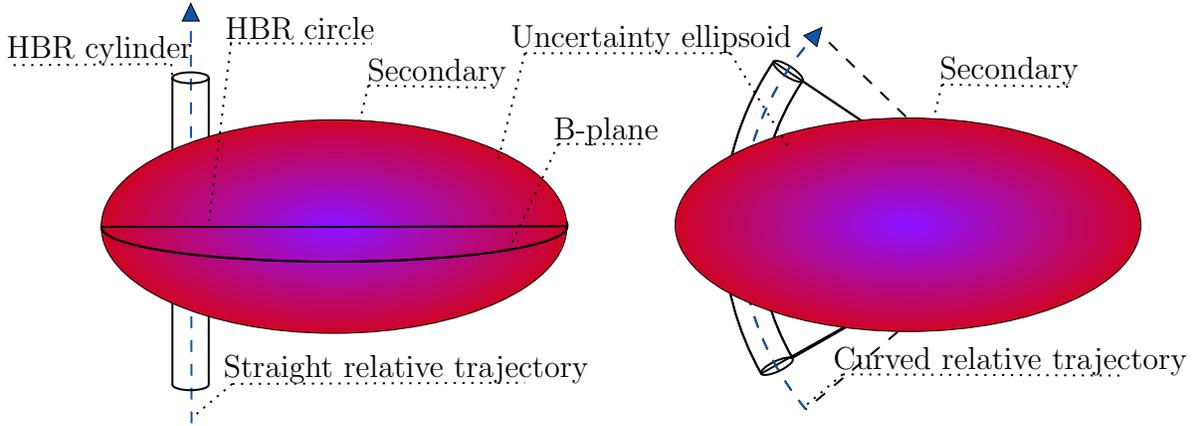}
    \caption{Graphical representation of short- (left) and long-term (right) encounters.}
    \label{fig:encounters}
\end{figure}

Two distinct encounter models have been explored in the existing literature: short- and long-term encounters. To distinguish them, the conjunction duration ($t_c$) is introduced by K. Chan in \cite{Chan2008Spacecraft}. It represents the time required for the primary object to cross the $1-\sigma$ relative position uncertainty ellipsoid in the $\eta$ B-plane direction at \ac{TCA} (which corresponds to the relative velocity direction) \cite{Chan_short}:
\begin{equation}
    t_c = \frac{2\sigma}{||\vec{v}_p -\vec{v}_s||}
\end{equation}
where $\vec{v}_p$ and $\vec{v}_s$ are the absolute velocities of the primary and secondary spacecraft, respectively.
By taking into account $t_c$ and the primary's orbital period $T_p$, a short-term encounter occurs when ${\epsilon=\displaystyle\frac{t_c}{T_p}\ll1}$ whilst a long-term one when ${\epsilon=\displaystyle\frac{t_c}{T_p}\approx 1}$. As sketched in \cref{fig:encounters}, in the former case, the relative motion inside the uncertainty ellipsoid is quasi-rectilinear, while it is curved in the latter.

The remainder of this chapter will present in further detail the short-term encounter techniques that are most commonly adopted for autonomous \ac{CAM} design and the \ac{IPoC} techniques used for long-term encounters. This is not intended as a review of currently available \ac{PoC} and {IPoC} approximations but as a necessary preliminary discussion of the models used in the \ac{CAM} methods presented in \cref{sec:NumericalCAM,sec:AnalyticalCAM}.
%%%%%%%%%%%%%%%%%%%%%%%%%%%%%%%%%
\subsubsection{Short-Term PoC}
\label{subsubsec:ShortPoC}

\begin{figure}[tb!]
    \centering
    \tikzset{every picture/.style={line width=0.75pt}} %set default line width to 0.75pt        

\begin{tikzpicture}[x=0.75pt,y=0.75pt,yscale=-1,xscale=1]
%uncomment if require: \path (0,465); %set diagram left start at 0, and has height of 465

%Straight Lines [id:da1719990401256153] 
\draw [line width=1.5]  [dash pattern={on 5.63pt off 4.5pt}]  (370,238) -- (371,121.67) ;
%Straight Lines [id:da6056213587618334] 
\draw [line width=1.5]  [dash pattern={on 5.63pt off 4.5pt}]  (371,121.67) -- (282.75,87.5) ;
%Straight Lines [id:da21160951945207707] 
\draw [color={rgb, 255:red, 208; green, 2; blue, 27 }  ,draw opacity=1 ][fill={rgb, 255:red, 208; green, 2; blue, 27 }  ,fill opacity=1 ][line width=0.75]    (283,203) -- (236.77,29.9) ;
\draw [shift={(236,27)}, rotate = 75.05] [fill={rgb, 255:red, 208; green, 2; blue, 27 }  ,fill opacity=1 ][line width=0.08]  [draw opacity=0] (8.93,-4.29) -- (0,0) -- (8.93,4.29) -- cycle    ;
%Straight Lines [id:da7666646074688507] 
\draw [color={rgb, 255:red, 74; green, 144; blue, 226 }  ,draw opacity=1 ][fill={rgb, 255:red, 208; green, 2; blue, 27 }  ,fill opacity=1 ][line width=1.5]    (283,203) -- (368.06,124.38) ;
\draw [shift={(371,121.67)}, rotate = 137.25] [fill={rgb, 255:red, 74; green, 144; blue, 226 }  ,fill opacity=1 ][line width=0.08]  [draw opacity=0] (11.61,-5.58) -- (0,0) -- (11.61,5.58) -- cycle    ;
%Straight Lines [id:da47921981090746657] 
\draw [line width=1.5]    (283,203) -- (408.29,253.5) ;
\draw [shift={(412,255)}, rotate = 201.95] [fill={rgb, 255:red, 0; green, 0; blue, 0 }  ][line width=0.08]  [draw opacity=0] (11.61,-5.58) -- (0,0) -- (11.61,5.58) -- cycle    ;
%Shape: Circle [id:dp7100792087747478] 
\draw  [fill={rgb, 255:red, 74; green, 144; blue, 226 }  ,fill opacity=1 ] (367.67,121.67) .. controls (367.67,119.83) and (369.16,118.33) .. (371,118.33) .. controls (372.84,118.33) and (374.33,119.83) .. (374.33,121.67) .. controls (374.33,123.51) and (372.84,125) .. (371,125) .. controls (369.16,125) and (367.67,123.51) .. (367.67,121.67) -- cycle ;
%Straight Lines [id:da2872937414839446] 
\draw [line width=1.5]    (283,203) -- (139.87,240) ;
\draw [shift={(136,241)}, rotate = 345.51] [fill={rgb, 255:red, 0; green, 0; blue, 0 }  ][line width=0.08]  [draw opacity=0] (11.61,-5.58) -- (0,0) -- (11.61,5.58) -- cycle    ;
%Straight Lines [id:da7614184023654675] 
\draw [line width=1.5]    (283,203) -- (283,33) ;
\draw [shift={(283,29)}, rotate = 90] [fill={rgb, 255:red, 0; green, 0; blue, 0 }  ][line width=0.08]  [draw opacity=0] (11.61,-5.58) -- (0,0) -- (11.61,5.58) -- cycle    ;
%Straight Lines [id:da5741149422048625] 
\draw [color={rgb, 255:red, 74; green, 144; blue, 226 }  ,draw opacity=1 ][fill={rgb, 255:red, 208; green, 2; blue, 27 }  ,fill opacity=1 ][line width=0.75]    (283,203) -- (149.98,52.25) ;
\draw [shift={(148,50)}, rotate = 48.58] [fill={rgb, 255:red, 74; green, 144; blue, 226 }  ,fill opacity=1 ][line width=0.08]  [draw opacity=0] (8.93,-4.29) -- (0,0) -- (8.93,4.29) -- cycle    ;
%Straight Lines [id:da6917854388377487] 
\draw [color={rgb, 255:red, 65; green, 117; blue, 5 }  ,draw opacity=1 ][fill={rgb, 255:red, 208; green, 2; blue, 27 }  ,fill opacity=1 ][line width=1.5]    (283,203) -- (198.87,224.99) ;
\draw [shift={(195,226)}, rotate = 345.35] [fill={rgb, 255:red, 65; green, 117; blue, 5 }  ,fill opacity=1 ][line width=0.08]  [draw opacity=0] (11.61,-5.58) -- (0,0) -- (11.61,5.58) -- cycle    ;
%Shape: Circle [id:dp8022912382599631] 
\draw  [fill={rgb, 255:red, 208; green, 2; blue, 27 }  ,fill opacity=1 ] (279.67,203) .. controls (279.67,201.16) and (281.16,199.67) .. (283,199.67) .. controls (284.84,199.67) and (286.33,201.16) .. (286.33,203) .. controls (286.33,204.84) and (284.84,206.33) .. (283,206.33) .. controls (281.16,206.33) and (279.67,204.84) .. (279.67,203) -- cycle ;
%Straight Lines [id:da21867553564150388] 
\draw [color={rgb, 255:red, 65; green, 117; blue, 5 }  ,draw opacity=1 ][fill={rgb, 255:red, 208; green, 2; blue, 27 }  ,fill opacity=1 ][line width=1.5]    (236,27) -- (151.87,48.99) ;
\draw [shift={(148,50)}, rotate = 345.35] [fill={rgb, 255:red, 65; green, 117; blue, 5 }  ,fill opacity=1 ][line width=0.08]  [draw opacity=0] (11.61,-5.58) -- (0,0) -- (11.61,5.58) -- cycle    ;

% Text Node
\draw (262.33,26) node [anchor=north west][inner sep=0.75pt]    {$\xi $};
% Text Node
\draw (383.33,250) node [anchor=north west][inner sep=0.75pt]    {$\zeta $};
% Text Node
\draw (151,240) node [anchor=north west][inner sep=0.75pt]    {$\eta $};
% Text Node
\draw (192.33,57) node [anchor=north west][inner sep=0.75pt]    {$\textcolor[rgb]{0.82,0.01,0.11}{v_s(t_CA)}$};
% Text Node
\draw (126.33,80) node [anchor=north west][inner sep=0.75pt]  [color={rgb, 255:red, 74; green, 144; blue, 226 }  ,opacity=1 ]  {$\textcolor[rgb]{0.29,0.56,0.89}{v_p(t_{CA})}$};
% Text Node
\draw (290,129) node [anchor=north west][inner sep=0.75pt]  [color={rgb, 255:red, 74; green, 144; blue, 226}  ,opacity=1 ]  {$\textcolor[rgb]{0.29,0.56,0.89}{r_{rel}(t_{CA})}$};
% Text Node
\draw (167.33,10) node [anchor=north west][inner sep=0.75pt]  [color={rgb, 255:red, 74; green, 144; blue, 226 }  ,opacity=1 ]  {$\textcolor[rgb]{0.25,0.46,0.02}{v_{rel}(t_{CA})}$};
% Text Node
\draw (249,213) node [anchor=north west][inner sep=0.75pt]   [align=left] {\textcolor[rgb]{0.82,0.01,0.11}{Secondary}};
% Text Node
\draw (361,92) node [anchor=north west][inner sep=0.75pt]   [align=left] {\textcolor[rgb]{0.29,0.56,0.89}{Primary}};
% Text Node
\draw (204,229) node [anchor=north west][inner sep=0.75pt]  [color={rgb, 255:red, 74; green, 144; blue, 226 }  ,opacity=1 ]  {$\textcolor[rgb]{0.25,0.46,0.02}{v_{rel}(t_{CA})}$};

\end{tikzpicture}
    \caption{B-plane construction (adapted from Ref. \cite{Armellin2021}). $v_{rel}$ and $r_{rel}$ are the relative velocity and the position of the primary in the B-plane, respectively.}
    \label{fig:bplane}
\end{figure}
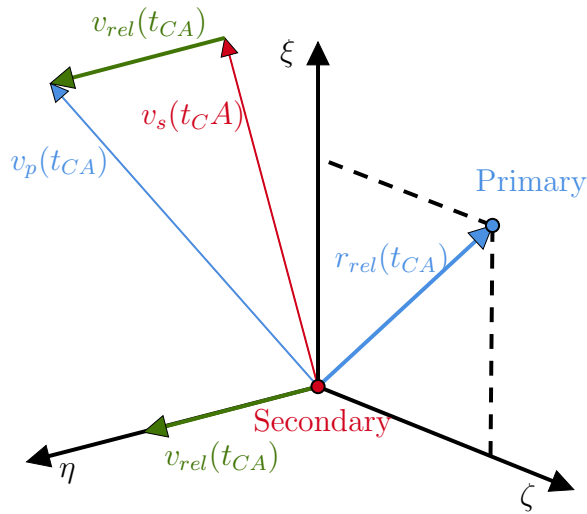

In this model, several assumptions can be adopted. Namely, (i) the positions of the two objects are considered uncorrelated Gaussian random variables, (ii) the encounter region is defined as $n$ standard deviations ($n\sigma$), (iii) the relative velocity is sufficiently high for the relative path to be rectilinear inside the combined uncertainty ellipsoid, (iv) no velocity uncertainty is involved, and (v) the objects are approximated as spherical, each with its specific radius. The two radii are commonly lumped together in the combined \ac{HBR}, which is the only parameter that retains information on the size.
The B-plane reference frame is a commonly adopted representation of short-term conjunction events. As illustrated in \cref{fig:bplane}, its centre is set on one of the involved objects (usually the secondary). 

The $\xi$ axis is aligned with the cross product of the velocities of the two spacecraft; the $\eta$ axis is directed along the relative velocity; finally, the $\zeta$ axis completes the right-handed triad.
Since the B-plane is defined only at \ac{TCA}, the relative position always belongs to the $\xi-\zeta$ plane. In fact, the closest approach condition enforces
\begin{equation}
    \vec{r}_{rel}\cdot\vec{v}_{rel} = 0.
\end{equation}
where $\vec{r}_{rel}, \vec{v}_{rel} \in\mathbb{R}^3$ are the relative position and velocity. Moreover, since position uncertainties are considered uncorrelated Gaussians, they can be summed and centred on the primary object without loss of generality. These simplifications allowed the development of several short-term \ac{PoC} models over the years. 

All the models aim at simplifying the solution of the \ac{PoC} integral for Gaussian distributions
\begin{equation}
    \pc = \iint_{\mathbb{C}}\frac{1}{\sqrt{(2 \pi)^{3}\mathrm{det}(\vec{P})}} \mathrm{exp}\left(-\frac{(\vec{r}-\vec{\mu}_r)^\mathrm{T} \vec{P}^{-1} (\vec{r}-\vec{\mu}_r)}{2}\right)\mathrm{d}\xi \mathrm{d}\zeta
    \label{eq:pdf}
\end{equation}
where $\vec{r}\in\mathbb{R}^2$ is drawn from the relative position random variable, ${\vec{\mu}_r\in\mathbb{R}^2}$ is the mean and $\vec{P} \in \mathbb{S}_{++}^2,$ where  $\mathbb{S}_{++}^n = \left\{ \mathbf{X} \in \mathbb{S}^n \mid \mathbf{X} \succ 0 \right\}$
 with $\mathbb{S}^{n} $ the set of symmetric matrices, is the covariance matrix of the relative position in the B-plane; $\mathbb{C}$ is the integration region identified by the combined \ac{HBR} $\in\mathbb{R_+}$.

Foster and H. S. Estes designed a strategy for circular cross-sections that relies on numerical techniques \cite{Foster1992}. This method has gained prominence and is widely used within \ac{NASA} for its computational accuracy.
Alfriend \textit{et al.} performed a study to assess the risk of collision for the \ac{ISS} in 1999 \cite{Alfriend2000}. The technique requires that the characteristic length of the relative position uncertainty exceeds the combined \ac{HBR} by at least one order of magnitude. By assuming a constant \ac{PDF} within the combined hard body sphere, the \ac{PoC} integral calculation is greatly simplified. Among the formulas in this section, \cite{Alfriend2000} stands out for its coarse approximation while effectively linking \ac{PoC} and \ac{SMD}. In the same work, Alfriend \textit{et al.} also estimate the maximum \ac{PoC} for an encounter, accounting for uncertainties in the covariance matrix, which provides an upper enveloping function for the actual \ac{PoC}, and serving as a conservative estimate.

The \ac{MECSA}, initially postulated by Chan in 2004 \cite{Chan2004Short,Chan2004International}, is an analytical expression that can be transposed to cross-sections of various shapes. Excelling in accuracy, the method is particularly effective for shapes such as circles, ellipses, rectangles, hexagons, octagons, and any cross-sections with bilateral symmetry. It can also be adapted to triangles and complex polygons—including concave shapes and those with internal voids—and extended to cross-sections bounded by curves.
Bai and Chen \cite{Bai2010Space} similarly derive a solution by transforming the anisotropic \ac{PDF} into an isotropic one and representing the resulting ellipsoidal hard body shape as circular, reducing the integration to a circular \ac{PDF} over a circular domain. Serra \textit{et al.} \cite{Serra2016Fast} further refine Chan's method by introducing an analogous to the \ac{IPoC} technique from \cite{Serra2015}. Leveraging the Laplace transform and properties of the D-infinite function, they develop a novel analytical integration strategy, with Chan's formulation emerging as a special case.
The ``Line" Integral method, introduced by R. P. Patera in 2001 \cite{Patera2001General} and subsequently modified in 2004, is a numerical approach that theoretically applies to cross-sections of any shape. However, it faces challenges with changes in viewing angles. Similarly, the Error Function Method, conceived by Alfano in 2005 \cite{Alfano2005Numerical}, works for any cross-section shape but also struggles with viewing angle variations and encounters computational difficulties when handling extremely low collision probabilities.

\subsubsection{Instantaneous PoC}
\label{subsubsec:InstantPoC}

Computing the \ac{PoC} for long-term encounters is computationally demanding, with velocity uncertainty further complicating \ac{CAM} planning. To overcome these challenges, \ac{CAM} methods often replace the \ac{PoC} constraint with \ac{IPoC}, an instantaneous snapshot of \ac{PoC} that provides a less rigorous risk representation. Academic interest in \ac{IPoC} has been limited, and fewer approximation methods exist compared to \ac{PoC}. However, from a practical standpoint, \ac{IPoC} remains a valuable tool for \ac{CAM} design. In Ref. \cite{Pavanello2024}, the authors demonstrate the validity of substituting \ac{PoC} with \ac{IPoC}: they show that, in typical scenarios, a sharp increase in the \ac{PoC} function is often accompanied by a rise in \ac{IPoC}.

While in most short-term encounter methods \ac{PoC} only depends on the mean state and covariance of the two objects at \ac{TCA}, \ac{IPoC} in long-term encounters is typically studied over an extended time window. The formal definition of \ac{IPoC} reads

\begin{equation}
\ipc(t) =\iiint_{\mathbb{S}} \frac{1}{(2 \pi)^{3 / 2}\mathrm{det}(\vec{P}(t))^{1/2}} \mathrm{exp}\left(-\frac{(\vec{r}-\vec{\mu}_r)^\mathrm{T} \vec{P}^{-1} (\vec{r}-\vec{\mu}_r)}{2}\right) \mathrm{d} V
\label{eq:ipceq}
\end{equation}
where $\ipc(\cdot):\mathbb{R^3}\rightarrow\mathbb{R_+}$ is \ac{IPoC} and the relative position random variable is defined in the tridimensional space; $\mathbb{S}$ is the spherical integration region identified by the combined \ac{HBR}. The covariance can be propagated using one of the methods illustrated in \cref{subsec:UP} to determine the time evolution of the integral. In the following discussion, the methods that are used in the works reviewed in \cref{subsec:NumericalCAM_Long} are briefly presented.

Chan proposed two methods to compute \ac{IPoC} either in an approximate analytical or in a semi-analytical manner \cite{Chan2008Spacecraft}. In the first case, the initial Gaussian distribution is transformed into a central chi-squared distribution, enabling \ac{IPoC} computation using an error function and chi-squared parameters, thus eliminating integration. This results in a fast method, though it can be imprecise in some instances. In the second method, the integration volume is switched from an ellipsoid into an equivalent sphere, and the integration is performed over a dimension analytically and over the second numerically.
In a similar fashion, Serra \textit{et al.} devised a method to define \ac{IPoC} via an infinite series of polynomial coefficients: the method is rapid and can achieve more accurate results than Chan's \cite{Li2022}.

Numerical methods are typically more taxing compared to Chan's and Serra's, but they can yield more accurate results. Notably, Zhang et al. \cite{Zhang2020} introduced a model which approaches the problem of the 3D integral by separating it into the multiplication of three-line integrals. To do so, three consecutive transformations are applied: (i) the ellipsoid shape of the covariance is contracted into a sphere, thus transforming the spherical hard body into an ellipsoid; (ii) the problem is rotated to express it in a reference frame where the covariance matrix is diagonal; (iii) the volume of the ellipsoidal shape of the hard body is transformed into a cuboid of equal volume. Integrating the uncorrelated \ac{PDF} in the three linear dimensions becomes a computationally efficient task.

Lastly, in \cite{Pavanello2024}, an approximation of \cref{eq:ipceq} that considers a constant value of the \ac{PDF} inside the combined hard body sphere is proposed. To compute this approximation, the same approximation from \cite{Alfriend2000} is used: only the \ac{PDF} value at the centre of the sphere is considered, which allows one to make the triple integral independent from \ac{SMD}. Moreover, using the same approach, an instantaneous analogous to Alfriend's maximum \ac{PoC} can be defined, accounting for the uncertainty in the covariance matrix estimation. The profile of maximum \ac{IPoC} provides an upper envelope to the time evolution of \ac{IPoC}, as shown in Ref. \cite{Pavanello2024}.

% This results in the simple formula
% \begin{equation}
%    \ipc(t) = \sqrt{\frac{2}{\pi \mathrm{det}(\vec{P}(t))}}\frac{R^3}{3}\mathrm{e}^{-\smd(t)/2}.
%    \label{eq:constipc}
% \end{equation}

% \begin{equation}
%     \ipc(t) = \frac{(\sqrt{2} R)^3 }{3\mathrm{e}^1 \smd(t) \sqrt{\pi \mathrm{det}\left(\vec{P}(t)\right)}}.
%     \label{eq:maxIpc}
% \end{equation}

%%%%%%%%%%%%%%%%%%%%%%%%%%%%%%%%%%%%%%%%%%%%%%%%%%%%%%%%%%%%%%%%%%%%%%%%%%%
% Numerical CAM
%%%%%%%%%%%%%%%%%%%%%%%%%%%%%%%%%%%%%%%%%%%%%%%%%%%%%%%%%%%%%%%%%%%%%%%%%%%
\section{Numerical CAM Design Methods}
\label{sec:NumericalCAM}

In this section, an overview of the \ac{CAM} methods that rely on numerical optimisation is given. In general, this includes all those methods that are based on the discretisation and successive numerical solution of the original \ac{CAM} \ac{OCP}. 
Commonly, direct methods are more flexible than indirect ones, and they allow for the inclusion of relevant operational and saturation constraints. They can be adapted to both impulsive and low-thrust dynamics and often include multiple safety constraints (e.g., \ac{MD}, \ac{SMD}, \ac{PoC}). However, this adaptability reduces computational efficiency, making them slower than indirect methods. Additionally, numerical optimisation algorithms, in general, have no convergence guarantees to the global optimum; therefore, their applicability in autonomous onboard applications should be carefully considered.

Literature contributions are grouped into three sections: in the first one, the methods that only treat short-term encounters are presented; in the second one, algorithms that can treat long-term encounters are included; in the third one, the discussion shifts to methods that address multiple short-term encounters.

%forse servono due tabelle? C'e' da considerare anche direzione di spinta, tempo computazionale, altri constraint operativi, multiple encounters, uncertainty propagation, ottimizzazione del tempo di manovra, fuel-loss, orbit (circular o eccentric) ...

%%%%%%%%%%%%%%%%%%%%%%%%%%%%%%%%%%%%%%%%%%%%%%%%%%%%%%%%%%%%%%%%%%
\subsection{Short-Term Encounters}
\label{subsec:NumericalCAM_Short}

The methods discussed in this section are reported briefly in \cref{tab:summaryShort}.

\begin{table}[tb!]
\footnotesize
    \centering
    \caption{Summary of numerical \ac{CAM} methods dealing with single short-term encounters}
    \label{tab:summaryShort}
    \begin{tabular}{lllcccl}
    \hline
        Ref.                                                & Approach          & \ac{PoC} model  & Low-thrust & Op. constraints & Free-time & Dynamics     \\ \hline
        \cite{Patera2003}                                   & Newton-Rhapson    & Chan            &  \xmark    & \xmark          & \xmark    & Keplerian    \\
        \cite{Pineiro2011}                                  & Pseudospectral    & \xmark          &  \cmark    & \xmark          & \cmark    & Keplerian    \\   
        \cite{Lee2012Collision}                             & Heuristic (GA)    & Not defined     &  \xmark    & \cmark          & \cmark    & Keplerian    \\
        \cite{Greco2021}                                    & Heuristic (SST)   & \xmark          &  \xmark    & \xmark          & \xmark    & Hi-fidelity  \\
        \cite{Armellin2021}                                 & SCP               & Alfriend        &  \xmark    & \xmark          & \cmark    & $J_2$        \\
        \cite{Dutta2022,Dutta2023a}                         & SOCP              & Alfriend        &  \cmark    & \cmark          & \cmark    & Hi-fidelity  \\
        \cite{SANCHEZ20232627}                              & minimax           & Chan            &  \cmark    & \xmark          & \xmark    & Keplerian    \\ 
        \cite{Pavanello2024Multiple,Pavanello2024LowThrust} & SCP               & Chan            &  \cmark    & \cmark          & \cmark    & Hi-fidelity  \\
        \cite{Pavanello2024Recursive}                       & PP                & Chan            &  \cmark    & \xmark          & \xmark    & $J_2$        \\
        \cite{Lopez2024}                                    & TFC               & Chan            &  \cmark    & \xmark          & \xmark    & Keplerian    \\
         \hline
    \end{tabular}
\end{table}

Patera and Peterson \cite{Patera2003} introduce an impulsive \ac{CAM} optimization method that assumes small manoeuvres. The firing direction is aligned with the \ac{PoC} gradient, while the magnitude is determined iteratively by expressing the closest approach displacement as a linear function of the $\Delta v$. A Newton-Rhapson algorithm is used to find the minimum of \ac{PoC} as a function of the impulse. This decoupling reduces the optimization to a one-dimensional problem.

In Ref. \cite{Pavanello2024Recursive}, the authors extend this approach by casting the problem as a \ac{PP}. They employ \ac{DA} to express \ac{PoC} as a function of the control input - whether low-thrust or impulsive - and across an arbitrary number of nodes, yielding a polynomial constraint. In addition, a quadratic objective function is used, effectively turning the problem into a \ac{PP}. Like in Patera and Peterson's method, in the first instance, the optimal control is approximated by firing in the direction of the gradient of \ac{PoC}, \textit{de facto}, eliminating the need for optimisation. Subsequently, the higher-order terms of the polynomial \ac{PoC} constraint are incorporated into the algorithm by linearising the tensors of high-order derivatives using the solution of the previous order. This allows one to successively change the magnitude and the direction of thrust until convergence.
This approach also accounts for the \ac{TCA} shift due to the manoeuvre, disregarded by most of the other methods. The achieved results are promising, with computation times in the order of tenths of a second that mainly depend on the expansion order and the number of impulses (or low-thrust arcs).
Ref. \cite{Pavanello2024IAC} extends the method to handle multiple constraints. Since it is used to address the multiple encounters scenario, its discussion follows in \cref{subsec:NumericalCAM_Multiple}.

Lee \textit{et al.} \cite{Lee2012Collision} are among the first to propose a \ac{GA} to find the optimal \ac{CAM}. They include \ac{SK} constraints for \ac{GEO} satellites and ground-track maintenance ones for \ac{LEO} and the option to design the manoeuvre a long time before \ac{TCA}. This approach has been proven unnecessarily conservative by research in the following decade, as most operators prefer to manoeuvre close to the conjunction. Their method is a precursor to following studies on multiple encounters; therefore, a more detailed discussion follows in \cref{subsec:NumericalCAM_Multiple}.
Sánchez and Vasile \cite{SANCHEZ20232627} implement a linear approximation of dynamics with impulsive manoeuvres. In their methodology, \ac{PoC} is computed under a short-term encounter assumption and integrates epistemic uncertainty in key parameters - such as \ac{MD} and covariance components. By accounting for this uncertainty, the encounter produces a range of \ac{PoC} values reflecting variations in the parameters. Thus, the problem reduces to determining the minimum manoeuvre needed to minimize the maximum collision risk—specifically under the worst-case scenario within the uncertain conjunction family.
The optimal manoeuvre direction is calculated for the worst-case \ac{PoC} arising from epistemic uncertainty along the direction of the minimum eigenvector of the combined covariance. After each manoeuvre computation, the conjunction geometry is updated, necessitating a recalculation of the extremal value of \ac{PoC} arising from the uncertain variable set. This iterative approach continues until the algorithm converges and the maximum \ac{PoC} stabilises. A line search strategy is integrated into the iteration to adjust the magnitude of the $\Delta v$, ensuring that the maximum \ac{PoC} post-manoeuvre falls below the selected risk mitigation threshold.
Additionally, the paper extends this framework to a heuristic, suboptimal strategy for low-thrust scenarios. This is achieved by discretising firing intervals and assuming a fixed optimal thrust direction for each interval, which is determined using the impulsive strategy for the respective sub-interval. The manoeuvre design strategy remains similar to the impulsive case, with the exception that the final position of the primary in the conjunction plane after the manoeuvre is determined via a semi-analytical propagation and the optimal direction within each arc is fixed. This fixed-direction assumption leads to the suboptimal nature of the method.
As with the impulsive case, the thrust arc magnitudes or manoeuvre amplitudes are adjusted via line search during the iterative solution of the minimax problem.

Pi{\~{n}}eiro \textit{et al.} \cite{Pineiro2011}  propose a low dimensional \ac{NLP} solution through a direct collocation method based on pseudo-spectral approximation of the controls and states. The solution aims at minimising the low-thrust fuel consumption while meeting a final state constraint linked to \ac{MD}. In fact, the approach avoids \ac{PoC} estimation and instead selects an appropriate semi-major axis variation to ensure safety against state and manoeuvre execution uncertainties at \ac{TCA}. The preliminary results focus solely on planar transfers from circular reference \ac{LEO} orbits, facilitating validation against classical Hohmann transfers.

A different take to the problem is proposed by Greco \textit{et al.} \cite{Greco2021} who adopt a Bayesian paradigm and importance sampling to determine the optimal impulsive \ac{CAM} in the presence of uncertainties in its execution (i.e., using Gates model). In particular, \ac{PoC} is computed as the expectation of the indicator function that assumes value $1$ if the \ac{MD} is smaller than a threshold and 0 otherwise. To obtain the optimal manoeuvre in this framework, the authors resort to a stochastic shooting approach with direct transcription.

Recently, convex optimisation has emerged as a promising approach for \ac{CAM} design. This branch of direct optimisation is particularly attractive for solving trajectory design problems due to its established success in other aerospace applications \cite{Liu2017,Malyuta2021Tutorial,Wang2024}. These studies highlight the advantages of formulating the convexified problem as a \ac{SOCP}, which facilitates the minimisation of the $l^2$ norm of the control history - a critical objective in spacecraft trajectory optimisation. To achieve this, all non-convex elements of the original \ac{OCP} must be convexified \cite{Boyd2}, including the fuel-optimal objective function, the non-linear dynamics constraint, and the \ac{PoC} constraint.

In Ref. \cite{Armellin2021}, the author applies the projection and linearisation method proposed by Mao in Ref. \cite{Mao2017} and later refined in Ref. \cite{Mao2021} to the short-term \ac{CAM} problem. This approach utilises \ac{DA} for automatic differentiation of the multi-impulsive dynamics and implements two nested iteration cycles to iteratively linearise the \ac{PoC} constraint based on a fixed dynamics reference solution. Furthermore, as in Ref. \cite{Pavanello2024Recursive}, the expansion of the conjunction time map is used to understand how the manoeuvre influences \ac{TCA}. The author also emphasises the significance of the initial trajectory guess, noting that it can cause the optimiser to converge to a local optimum on the avoidance ellipse. Consequently, while the method reliably produces locally optimal solutions, it does not guarantee a globally optimal solution to the original \ac{OCP}. Extensive simulations show that accurate results are achieved with sub-second performance, with most cases converging in just two \ac{SCP} iterations. As with all \ac{SCP} methods, there is no theoretical proof for the convergence of the algorithm. However, by means of a vast simulation campaign, the author shows that convergence is reached in most of the considered cases. Armellin's method is the basis for the development of a low-thrust-based \ac{SCP} algorithm that is described in \cref{subsec:NumericalCAM_Multiple}, but can, as a subcase, treat single encounters \cite{Pavanello2024Multiple}. The same method is adapted in Refs.\cite{Pavanello2024LowThrust,Pavanello2024LT} to address conjunctions that involve a continuously propelled low-thrust spacecraft. In this case, multiple objective functions are considered, which have the general aim of minimizing the deviation from the nominal propelled trajectory. A switch-off strategy and two minimum steering ones are proposed. While the first one is the simplest to implement and actuate, the minimum deviation strategies show good results in avoiding collision and re-routing the spacecraft to its nominal trajectory over a period of time of 10 to 40 hours, using non-demanding attitude changes.

%%%%%%%%%%%%%%%%%%%%%%%%%%%%%%%%%%%%%%%%%%%%%%%%%%%%%%%%%%%%%%%%%%%%%%%%%

Although all of the previously mentioned works are based on the use of some form of \ac{PoC}, Dutta and Misra \cite{Dutta2022} solve the short-term \ac{CAM} using the \ac{IPoC}. By employing the 3D \ac{COLA} ellipsoid, they use the same approximation of \cite{Alfriend2000} to find a constant-\ac{PDF} approximation of \ac{IPoC}.
This allows them to define the \ac{CAM} constraint on \ac{SMD}, similarly to the conversion from \ac{PoC} in Ref. \cite{Armellin2021}. Differently from the previously cited method, however, Dutta and Misra eliminate the need for successive projection and linearisation of the \ac{CAM} constraint by expanding the \ac{SMD} around a reference point. This linearisation leads to a more conservative solution, so in general, the method always finds a sub-optimal solution. However, given that the algorithm encompasses a single \ac{SOCP} instance without sequences of convexifications, it is guaranteed to converge to the local optimum of the original \ac{OCP}.
Nonetheless, the use of \ac{IPoC} instead of \ac{PoC} is a limiting factor, and the reduction of the $\Delta v$ is obtained at the price of a lower reduction of the risk.
% Nonetheless, the operational feasibility of the approach is undermined by another aspect, namely, the use of \ac{IPoC} instead of \ac{PoC}: not using the short-term B-plane model (i.e., not considering the ellipsoidal infinite cylinder that derives from the short-term assumption) does not take into account the propagation of the objects in the vicinity of \ac{TCA}. Therefore, the claim of achieving a lower $\dv$ may be invalidated as the authors are not genuinely reducing \ac{PoC}, but rather its less conservative instantaneous counterpart.
In Ref. \cite{Dutta2023a}, the authors extend the formulation to include a non-linear propagation of the covariance of the two objects. The covariance evolution is modelled using a Monte Carlo approach for the secondary and a \ac{UT} for the primary. Given the non-Gaussian nature of the uncertainty cloud at \ac{TCA}, \ac{SMD} cannot be used, and an analogous metric based on the cumulative density function is adopted to ensure that the sigma points of the \ac{UT} respect the avoidance requirement.
Because of the high computational load required by the Monte Carlo and by the \ac{UT}, the method is expensive, with runtimes that reach $2$~\si{hours}.

At last, Ref. \cite{Lopez2024} explores the \ac{TFC} to design energy-optimal \ac{CAM}s, framing the problem as an unconstrained optimization. A comparison is made between the \ac{TFC}-based approach and the analytical approach from Ref. \cite{DeVittori2022}, evaluating performance in terms of computational cost, accuracy, and manoeuvre cost while also analysing the impact of various \ac{TFC} configuration parameters. Despite \cite{DeVittori2022}, detailed in \cref{subsec:AnalyticalCAM_LowThrust}, offering higher computational efficiency, \ac{TFC} paves the way for incorporating any constraint without the need for dynamics and constraint convexification. 
% \textbf{CAM SHORT TERM CONFERENCE PAPER, TFC di David}

%%%%%%%%%%%%%%%%%%%%%%%%%%%%%%%%%%%%%%%%%%%%%%%%%%%%%%%%%%%%%%%%%%
\subsection{Long-Term Encounters}
\label{subsec:NumericalCAM_Long}

\begin{table}[tb!]
\footnotesize
    \centering
    \caption{Summary of numerical \ac{CAM} methods dealing with long-term encounters}
    \label{tab:summaryLong}
    \begin{tabular}{lllcccl}
    \hline
        Ref.                                        & Approach       & \ac{IPoC} model      & Low-thrust  & Op. constraints       & Free-time  &  Dynamics    \\ \hline
        \cite{Mueller2009}                          & LP             & \xmark               & \cmark      & \xmark   & \xmark     &  Keplerian   \\
        \cite{Mueller2013}                          & LP             & \xmark               & \cmark      & \cmark   & \cmark     &  $J_2$   \\
        \cite{Serra2015}                            & LP             & Serra                & \xmark      & \xmark   & \xmark     &  Keplerian   \\
        \cite{Pavanello2024,Pavanello2024Multiple}  & SCP + NLP      & Alfriend/Zhang       & \cmark      & \cmark   & \cmark     &  Hi-fidelity \\
        \cite{DeVittori2025}                        & SCP            & Alfriend             & \cmark      & \cmark   & \cmark     &  Hi-fidelity \\
 \hline
    \end{tabular}
\end{table}

Few attempts over the years have been conducted to devise \ac{CAM} strategies for long-term encounters, as reported in \cref{tab:summaryLong}. The first proposals are formulated in terms of \ac{LP} \cite{Mueller2009,Mueller2013} and deterministic disjunctive \ac{LP} \cite{Serra2015}. Mueller \cite{Mueller2009} develops a \ac{COLA} strategy that focuses on maintaining a minimum separation distance from the secondary spacecraft. This strategy involves approximating the nonconvex \ac{COLA} constraint through a linearisation around a point on the ellipsoid’s surface, specifically chosen on the semi-major axis, to ensure consistent along-track separation over time. Additionally, the amplitude of the oscillation of this separation is constrained, turning the problem into a \ac{COLA}-constrained relative navigation challenge. To simplify the relative motion, Mueller employs a linear time-varying system, transforming the \ac{CAM} problem into a \ac{LP} problem. To maintain linearity, the $l^1$ norm of the control is minimised, although this can lead to higher propellant consumption compared to the $l^2$ norm. Importantly, Mueller’s approach does not consider state uncertainty, as the \ac{KOZ} region is arbitrarily defined by the author. In Ref. \cite{Mueller2013}, the authors improved the algorithm by introducing a \ac{SK} linear constraint. The \ac{SK} requirement for \ac{GEO} is formulated in \ac{ECI} coordinates as a soft constraint. By varying the selected penalty weight in the objective function, the authors show how their approach can be flexible and operationally useful.

The first long-term encounter \ac{CAM} algorithm to consider state uncertainty, and, consequently, the covariance evolution, is the one proposed by Serra \textit{et al.} \cite{Serra2015} in 2015. They formulate a deterministic disjunctive \ac{LP} based on a linearization of the dynamics using the Yamanaka-Ankersen state transition matrix, which is also employed to propagate the covariance in a linear fashion. In this framework, since a relative motion model is used, the combined covariance is evaluated from \ac{TCA} to the desired final time. However, this may introduce inaccuracies linked to differences between the orbits of the two spacecraft, which lead to distinct state transition matrices. Due to the high computational load of a more accurate approach, \ac{PoC} is replaced with \ac{IPoC}, though this variable is not directly controlled.
Instead, a threshold value of \ac{IPoC} establishes a \ac{KOZ}. The \ac{KOZ} region is approximated as a convex polyhedron with an arbitrary number of faces. The deterministic disjunctive \ac{LP} results from the combination of the \ac{CAM} constraints on the polyhedron faces and the use of the $l^1$ norm of the control in the objective function. The \ac{IPoC} evolution of the manoeuvred trajectory is validated through the approximation defined by Ref. \cite{Serra2015}, finding that the manoeuvre always reduces the risk below an acceptable level.

In Refs. \cite{Pavanello2023Long,Pavanello2024}, the authors provide a more rigorous approach to the design of the manoeuvre than the one proposed in Ref. \cite{Serra2015}: they employ a \ac{SCP} formulation, wherein a succession of \ac{SOCP}s are solved iteratively until convergence. \ac{DA} is used to automatically linearise high-fidelity dynamics, and the state and covariances of the two spacecraft are propagated independently. Moreover, the covariance history of the primary is updated within each iteration, allowing for the refinement of the \ac{COLA} region. In the same work, the convexification approach for the \ac{PoC} constraint proposed in Ref. \cite{Armellin2021} for short-term encounters and used in Ref. \cite{Pavanello2024LT} is extended to long-term encounters. The \ac{IPoC} constraint is turned into a \ac{KOZ} one, which identifies a time-dependent iso-probability 3D ellipsoid. The \ac{KOZ} constraint is then iteratively approximated into a succession of hemispaces. Differently to short-term encounters \cite{Armellin2021,Pavanello2024LT}, the \ac{KOZ} constraint is applied across the entire window of interest, rather than just at the nominal \ac{TCA}. A linear constraint - imposed on the gradient of \ac{SMD} - is introduced to limit the sensitivity of \ac{IPoC} to state uncertainties, improving robustness against thrust misalignment or navigation and control errors.
\ac{SK} constraints are included both for \ac{LEO} and for \ac{GEO} scenarios during the low-thrust manoeuvre and at the end of it. In the \ac{LEO} case, a constraint on the final state forces the spacecraft to return to its nominal orbit; in \ac{GEO}, a linearised keep-in-box constraint bounds the spacecraft to respect a latitude-longitude requirement during the manoeuvre, and the optimal state that grants a long resident time inside the box is targeted at the end of the manoeuvre. This allows the operator to combine \ac{CAM} and \ac{SK} manoeuvres. In a following work, the method is adapted to deal with a non-linear evolution of the covariance of the secondary via the introduction of \ac{GMM} \cite{Pavanello2024Multiple}. When the hypothesis of Gaussianity is lost, i.e., if the propagation window is too vast and/or the initial uncertainty is high, this improvement allows for a more reliable representation of the \ac{KOZ}. Indeed, in this case, the \ac{KOZ} is not a single ellipsoid but a more complex shape defined by the union of multiple ellipsoids. In Ref. \cite{DeVittori2025}, the method is ulteriorly improved via introducing stochastic \ac{SK} policies tailored to longitude and latitude, ensuring adherence to \ac{GEO} slots with probabilistic guarantees. Two types of chance constraints are proposed, either neglecting or including longitude-latitude correlation.

%%%%%%%%%%%%%%%%%%%%%%%%%%%%%%%%%%%%%%%%%%%%%%%%%%%%%%%%%%%%%%%%%%
\subsection{Multiple Encounters}
\label{subsec:NumericalCAM_Multiple}
% Kim, Seong, Sanchez, Masson, Zeno

\begin{table*}[b!]
\footnotesize
    \centering
    \caption{Summary of numerical \ac{CAM} methods dealing with multiple encounters}
    \label{tab:summaryMultiple}
    \begin{tabular}{lllcccl}
    \hline
        Ref.                            & Approach           & \ac{PoC} model     & Low-thrust & Op. constraints       & Free-time  & Dynamics \\ \hline
        \cite{Kim2012}                  & Heuristic (GA)     & Alfano             &  \xmark    & \xmark                & \cmark     & Hi-fidelity \\
        \cite{Morselli2014Collision}    & Multi-objective    & Alfano             &  \xmark    & \cmark                & \cmark     & Keplerian \\
        \cite{Seong2015,Seong2016}      & Heuristic          & \xmark             &  \xmark    & \cmark                & \cmark     & Hi-fidelity \\
        \cite{Sanchez2021}              & min-max            & Full integration   &  \xmark    & \cmark                & \cmark     & Keplerian \\
        \cite{Masson2023}               & PP                 & Serra              &  \xmark    & \cmark                & \xmark     & Keplerian \\
        \cite{Bourriez2023}             & RL                 & Monte Carlo        &  \xmark    & \cmark                & \xmark     & Hi-fidelity  \\ 
        \cite{Pavanello2024Multiple}    & SCP                & Chan               &  \cmark    & \cmark                & \cmark     & Hi-fidelity \\
        \cite{Pavanello2024IAC}         & PP                 & Chan               &  \cmark    & \cmark                & \xmark     & $J_2$  \\ 
        \cite{Mu2024}                   & RL                 & Bai                & \xmark     & \xmark                & \cmark     & $J_2$ \\
        \hline
    \end{tabular}
\end{table*}

While most existing research focuses on single encounters, numerical methods are increasingly being used to study the \ac{CAM} optimisation problem in scenarios involving multiple encounters, as reported in \cref{tab:summaryMultiple}. The first of these attempts came in the early 2010s by Kim \textit{et al.} \cite{Kim2012} expanding on their work on single encounters \cite{Lee2012Collision}. These approaches do not rely on a rigorous mathematical formulation of an optimal control problem in a classical sense but rather use heuristic techniques to optimise a given objective function. The authors develop a \ac{GA} approach to deal with up to four conjunctions. The method finds a single impulsive manoeuvre at the optimal time to avoid consecutive conjunctions. It can either optimise the pure tangential or tridimensional impulse in \ac{RTN} coordinates. In the former case, though, the method cannot find a suitable solution to the problem. The \ac{COLA} requirement can be stated in terms of \ac{PoC} (computed with Alfano's formula \cite{Alfano2005Numerical}) or \ac{MD}. 

Seong and Kim \cite{Seong2015} improve on the previous method and explore different heuristic approaches to identify the most suitable one for the multiple encounters problem. They employ a high-fidelity dynamics model and make \ac{CAM} decisions based on a two-step assessment: a coarse evaluation using \ac{TLE} propagation, followed by a precise orbit determination. The thrust dynamics are modelled as impulsive, with \ac{MD} as the sole avoidance metric. The proposed approaches include a \ac{GA}, similar to Ref. \cite{Kim2012}, along with a particle swarm optimisation algorithm, a differential evolution algorithm, and a simulated annealing algorithm. 
These are different variants of heuristic approaches, where an initial population undergoes iterative refinement until converging to a sub-optimal final solution. In fact, no guarantee can be given to find an optimal solution of the \ac{CAM} \ac{OCP}. The efficacy of the approaches is demonstrated in a test case similar to \cite{Kim2012}, involving four consecutive conjunctions, wherein each method yields a control policy that raises \ac{MD} above the required threshold. Notably, the differential evolution algorithm delivers the highest fuel efficiency, while the \ac{GA} produces the least efficient one.

In a last research effort \cite{Seong2016}, Seong and Kim improve the \ac{GA} proposed in their previous work by introducing a second objective function, converting the problem into a multi-objective optimisation one. This allows them to combine the \ac{CAM} requirement with the ground-track maintenance one, typical for \ac{LEO} satellites. They propose 3 burn strategies: the first consisting of a single tangential burn, the second one of a single burn with free \ac{RTN} direction, and the third one with two burns with free \ac{RTN} direction. They find the third strategy to be the best in respecting the ground track evolution required for fourteen days. 
Since no comprehensive information is provided on \ac{TPoC}, these three methods typically underestimate the risk of collision. Moreover, \ac{GA}s (and heuristic methods in general) are dependent on the specific definition of the fitness function, which can be very problem-specific, and they are notoriously slower to converge compared to classic optimisation algorithms, with run-times between $20$ and $40$ \si{min}.

A different approach that led to a similar solution methodology is depicted in Ref. \cite{Morselli2014Collision}. In this work, the authors set out to determine an impulsive \ac{CAM} while looking not only at the nearest conjunction but also at all subsequent conjunctions in a predetermined time window (e.g., 7 days in \ac{LEO}). The authors leverage the \ac{SGP4} model and \ac{DA} to determine all subsequent closest approaches and linearly propagate the covariances at those epochs. Subsequently, \ac{PoC} calculated using Alfano's method \cite{Alfano2005Collision}, \ac{MD}, and fuel consumption are all considered in a multi-objective optimisation executed with a particle swarm technique. As a result, a Pareto optimal front of manoeuvres is obtained rather than a single optimal solution.

In Ref. \cite{Sanchez2021}, Sánchez and Vasile propose two methods for managing multiple encounters: one that derives a single manoeuvre based solely on the worst-case scenario, and another that treats each conjunction as an independent event. The first approach may yield a suboptimal impulse, whereas the second computes individual manoeuvres for each conjunction, requiring \ac{PoC} recalculations for subsequent events after each \ac{CAM}. This method primarily focuses on proposing and ranking various potential \ac{CAM}s under different operational and execution constraints. The authors employ rigorous multi-criteria decision-making techniques to evaluate collision risk and prioritise the candidate \ac{CAM}s. While this approach is advantageous for ground operators in selecting the most appropriate \ac{CAM} for a given scenario, its high computational demands and reliance on human intervention render it unsuitable for autonomous manoeuvres.

The first contribution in terms of pure direct optimisation on the topic of multiple conjunctions is brought by Massòn \textit{et al.} \cite{Masson2023}. They cast the problem as a quadratically constrained \ac{LP}, and they solve it by pre-allocating the manoeuvring opportunities. The dynamics are linearised using the Yamanaka-Ankersen state transition matrix, which allows them to directly relate the constraint in terms of \ac{MD}, \ac{SMD}, or \ac{PoC} to the impulsive manoeuvre opportunities. The optimisation problem is thus built with the $\dv$ history as optimisation variables: the $l^1$ norm of the impulses history is used as the objective function, the \ac{COLA} constraint is formulated as quadratic, independent from the specific collision metric utilised, and a \ac{SK} linear constraint is also included in terms of an \ac{RTN} box. In conclusion, the resulting quadratically constrained \ac{LP} can be seen as a \ac{PP}, which is solved using one of three methods. The first method employs Lasserre's hierarchy - a series of \ac{SDP} relaxations of the \ac{PP} - thereby guaranteeing a globally optimal solution. The second one utilises a semidefinite relaxation solved via randomisation, which is only capable of identifying local optima. The third approach relies on two branch-and-bound algorithms (IBEX and an ellipsoid-based variant). Once a history of optimal impulses has been established, an \textit{a posteriori} spread is applied to convert these impulses into a series of thrust arcs, enabling operators to realistically assess the manoeuvre. In a scenario with two consecutive conjunctions, all four methods produce nearly identical optimal impulses, with execution times ranging from 1 to 14 seconds. When dealing with more than two conjunctions, IBEX significantly outperforms the others, completing the task in less than one second. Additionally, Lasserre's hierarchy is relatively swift for fewer than eight conjunctions, with run times under 18~\si{s}.

Ref. \cite{Pavanello2024Multiple} presents a comprehensive approach to the multiple encounter problem, which uses a convex optimisation method similar to Ref. \cite{Pavanello2024}. The solution methodology is tailored to address the intricacies of multiple encounters, with consideration on \ac{TPoC} -- a novel concept introduced by \ac{NASA} \cite{Frigm2020} to quantify the contribution of multiple conjunctions to the overall \ac{PoC} -- and non-linear uncertainty propagation. The authors propose two alternative algorithms to comply with the \ac{TPoC} constraint. In the first one, the convex problem includes a linearised constraint on \ac{TPoC}, which is built automatically thanks to a first-order \ac{DA} expansion. The second method implements a sequence of \ac{SCP} interspersed with \ac{NLP}s to adaptively distribute probability limits across conjunctions, enforcing a linearised \ac{SMD} constraint to each. While the former method works best for a low number of conjunctions, becoming imprecise when more than four are present, the latter can treat any arbitrary number of conjunctions. The pipeline can also address repeating conjunctions with longer propagation times thanks to their generalised formalism. In this instance, the assumption of Gaussianity for the distributions is generally invalid, prompting the use of a \ac{GMM} to propagate the covariances with enhanced fidelity. Notably, since this work stems from Ref. \cite{Pavanello2024}, it can treat both short- and long-term multiple encounters. Numerical results demonstrate that for short-term encounters, outcomes with relative \ac{TPoC} accuracies below $10$\% can be achieved within acceptable time frames, ranging from $7$ to $30$~\si{s}, depending on the number of encounters.
When studying long-term encounters, instead, it is proven that increasing the number of kernels in the \ac{GMM} used to propagate uncertainties allows for a better representation of the evolution of the \ac{KOZ}, leading to more conservative manoeuvres.

In Ref. \cite{Pavanello2024IAC}, the authors extend the \ac{PP} approach from Ref. \cite{Pavanello2024Recursive} (described in  \cref{subsec:NumericalCAM_Short}) to handle multiple consecutive conjunctions. The core principle remains unchanged: progressively increasing the polynomial approximation order of the constraints to incorporate higher-order derivative terms. Unlike its initial version, the algorithm also handles inequality constraints, employing a state-of-the-art solver for each instance of the successive \ac{PP} formulations with linearised constraints. The authors show the method applied to multiple conjunctions, both using a comprehensive \ac{TPoC} approach and an individualistic \ac{PoC} or \ac{MD} one. Moreover, applications where a return to the nominal orbit constraint and a \ac{SK} constraint are showcased, showing promising results in terms of accuracy (relative constraint violations below $10^{-3}$) and efficiency (computation times below $0.5$~\si{s}). Nonetheless, as in its previous version, no theoretical proofs of convergence are given for the method, even though the simulations show a $99\%$ convergence rate over more than 2,000 test cases.

A notable mention must be reserved for those innovative approaches that employ AI techniques for \ac{CAM} design, such as Refs.\cite{Mu2024,Bourriez2023}. In these works, the authors adopt different reinforcement learning algorithms (i.e., proximal policy optimisation \cite{Mu2024} and Deep Recurrent Q-Network \cite{Bourriez2023}) to determine the best sequence of impulsive manoeuvres across a given planning horizon. In Ref. \cite{Mu2024}, the manoeuvre policy is determined in an infinite search space by maximising a specific reward composed of \ac{PoC}, energy consumption, and elapsed time (to encourage longer runs without failures); moreover, a sparse reward for reaching the final state while being compliant with constraints (e.g., \ac{MD}) is included.
Conversely, in Ref. \cite{Bourriez2023}, the possible impulsive manoeuvres are selected between a set of five predefined values (per direction) at each time, and a \ac{SK}-like constraint is embedded in the reward function through a penalty on the deviation from the reference trajectory.
In these preliminary studies, the debris population is synthetically generated to ensure multiple encounters occur during the \ac{CAM} planning window. Although this scenario is not realistic, it is essential to justify the application of reinforcement learning. Moreover, the state of the agent includes the orbital states of all debris, and in Ref. \cite{Mu2024}, the primary and secondary covariances are assumed to be constant and not propagated during the manoeuvre planning arc. Due to such strong assumptions, it appears clear that the level of maturity of this technology is insufficient, and the path that leads to the application of such a paradigm in real-world scenarios remains unclear.

%%%%%%%%%%%%%%%%%%%%%%%%%%%%%%%%%%%%%%%%%%%%%%%%%%%%%%%%%%%%%%%%%%%%%%%%%%%
% Analytical and Semi-Analytical Optimisation
%%%%%%%%%%%%%%%%%%%%%%%%%%%%%%%%%%%%%%%%%%%%%%%%%%%%%%%%%%%%%%%%%%%%%%%%%%%
\section{Analytical and Semi-Analytical CAM Design Methods}
\label{sec:AnalyticalCAM}

This section introduces analytical and semi-analytical \ac{COLA} algorithms, encompassing both impulsive and low-thrust propulsion strategies. These methods utilise closed-form solutions for \ac{CAM} problems, offering significantly improved computational efficiency compared to traditional numerical techniques.

The following discussion is organised into two parts: the first addresses impulsive strategies, while the second explores low-thrust manoeuvrers. Both analytical and semi-analytical formulations are examined within each subsection.

%%%%%%%%%%%%%%%%%%%%%%%%%%%%%%%%%
\subsection{Impulsive Methods}
\label{subsec:AnalyticalCAM_Impulsive}

Impulsive trajectory optimisation, typically used to design manoeuvrers for chemical thrusters, assumes that manoeuvres can be executed instantaneously with a change in velocity at discrete points in time. This kind of \ac{CAM}s often involves solving a sequence of \ac{BVP}, where the goal is to determine the optimal times, magnitudes, and directions of impulses to achieve a desired mission objective. In \cref{tab:AnalyticalSummary_Impulsive}, the methods discussed in this section are summarized.

\begin{table}[b!]
    \centering
    \caption{Summary of analytical and semi-analytical impulsive \ac{CAM} methods}
    \label{tab:AnalyticalSummary_Impulsive}
    \begin{tabular}{llllll}
    \hline
        Ref.                      & Approach        & Constraint & Dynamics    \\ \hline
        \cite{Bombardelli2014}         & Analytical      & MD         & Keplerian   \\
        \cite{Bombardelli2015}         & Analytical      & MD, Chan   & Keplerian   \\
        \cite{Gonzalo2021}             & Analytical      & Chan       & Keplerian   \\
        \cite{Abay2017}                & Semi-Analytical & Chan       & Hi-fidelity \\
        \cite{DeVittori2024capitolo6}  & Analytical      & Chan, MD   & Keplerian   \\
        \cite{RubioAntón2024}          & Analytical      & SMD, MD    & Keplerian   \\
        \hline
    \end{tabular}
\end{table}

The literature on this topic has a strong foundation with a key contribution by Bombardelli, who assessed the challenges of impulsive \ac{CAM}s assuming elliptical Keplerian orbits in \cite{Bombardelli2014}. Leveraging the B-plane dynamics at \ac{TCA}, closed-form analytical expressions are developed to increase \ac{MD} with a three-dimensional manoeuvre. The accuracy of these algorithms is tested under various orbital configurations, assessing its relevance for practical scenarios like space debris avoidance and satellite operations, where last-minute, non-tangential manoeuvres might be necessary.
Building upon this work, Bombardelli and Hernando-Ayuso in \cite{Bombardelli2015} present a high-accuracy analytical formulation to compute the \ac{MD} and \ac{PoC} of two approaching objects after an impulsive \ac{CAM}. The manoeuvre optimisation is turned into an eigenvalue problem coupled with a non-linear algebraic equation by exploiting the linear relation between applied impulse and objects' relative motion in the B-plane. As an optimisation variable, the cost of the manoeuvre is chosen in terms of $\Delta v$ magnitude to match imposed thresholds on either \ac{MD} or \ac{PoC}.

In Ref. \cite{Gonzalo2021}, Gonzalo \textit{et al.} introduce an analytical formulation for spacecraft \ac{CAM} with space debris, offering solutions for both maximum deviation and minimum \ac{PoC} scenarios. Employing Gauss' planetary equations and relative motion dynamics, the model maps manoeuvres at specific times into displacements at predicted close approaches. The formulation is extended to propagate covariance matrices, reducing the optimisation problem to an eigenvalue problem. Comparisons between \ac{MD} and \ac{PoC} methodologies are presented using data from the ESA's MASTER-2009 model. The analysis in the B-plane shows that both uncertainties and maximum deviation increase with lead time, limiting \ac{PoC} reduction over extended time frames.

A further advancement \cite{DeVittori2024capitolo6} extends impulsive \ac{CAM} strategies by introducing frameworks that seamlessly transition between \ac{ECI} or B-plane coordinates while incorporating constraints such as \ac{MD}, \ac{SMD} or \ac{PoC}. Single-impulse formulations are developed using \ac{NLP} to obtain closed-form solutions to the problem. These methodologies are further refined for the direct impact scenario, taking into account both fixed and free thrust direction constraints. A tri-impulsive analytical formulation is then introduced to return to the nominal trajectory: the first is applied at a predefined epoch to achieve lower \ac{PoC}, the second at \ac{TCA} to redirect the trajectory back toward the original orbit, and the third one at the point of intersection with the original orbit to achieve the \ac{SK} condition. These methods are tested in \ac{LEO}, \ac{MEO}, and \ac{GEO} regimes, demonstrating comparable computational efficiency and accuracy. Finally, R. Antón in \cite{RubioAntón2024} introduces a new formulation for designing impulsive Shared CAMs (ShCAM) for active vs. active conjunctions, utilising the Chan \ac{PoC} model and Bombardelli impulsive analytical formulation \cite{Bombardelli2014} as foundation for the manoeuvrer design and optimisation.

Concerning semi-analytical methods, R. Abay \cite{Abay2017} advances an optimal impulsive manoeuvre using \ac{SGP4} and \ac{TLE} data. The author analyses the relative dynamics of two colliding objects, assuming an instantaneous encounter, and tests various collision geometries. In all cases, the relative velocity and the \ac{MD} vectors are parallel, resulting in a worst-case scenario as the encounter occurs along the direction of the relative velocity vector. The method is tested and validated using SpaceTrack's catalogue to generate ephemerides of distinctive satellites. Successively, Ref. \cite{DeVittori2024capitolo6} elaborates a semi-analytical solution for the single-impulse \ac{CAM} leveraging the Dual Lagrangian function, followed by three semi-analytical methods for the two-impulses cases with return to the original orbit in the two-body problem. A first technique exploits the post-manoeuvre trajectory, letting the spacecraft coast the new orbit after \ac{TCA} and then return to the original one by firing again at the starting point (intersection between original and new orbit). The second method employs a phasing orbit, where the user defines the increase or decrease in the orbital period after the manoeuvre, and then an optimisation can be carried out to reduce the total $\Delta v$ while matching the prescribed safety metric. The last technique constructs a fictitious orbit by varying some Keplerian elements of the original one and then exploits the state transition matrix to find an orbit that has initial and final points on the original one. The methods are tested in \ac{LEO}, \ac{MEO}, and \ac{GEO} regimes, showing comparable computational times and accuracy between different orbital regimes.

%%%%%%%%%%%%%%%%%%%%%%%%%%%%%%%%%
\subsection{Low-Thrust Methods}
\label{subsec:AnalyticalCAM_LowThrust}

Low-thrust trajectory optimisation involves continuous thrusting over long durations, typically using electromagnetic propulsion systems, solar sails, or similar mechanisms. A key feature of low-thrust propulsion is the gradual but continuous acceleration, allowing for extended mission lifespans due to higher specific impulse and available $\Delta v$. This optimisation is particularly suited for \ac{CAM}s with sufficient warning time or for integrating them with \ac{SK} routines during close encounters. A noteworthy distinction exists between \ac{EO} and \ac{FO} formulations. The former minimises energy usage without thrust limitations, while the latter typically follows a bang-bang profile, maintaining constant acceleration throughout the manoeuvring windows. Moreover, \ac{EO} strategies often serve as a starting point for \ac{FO} optimisations through techniques like homotopy transformations \cite{Wang2023}. In \cref{tab:AnalyticalSummary_LowThrust}, the methods discussed in this section are summarized.

\begin{table}[H]
    \centering
    \caption{Summary of analytical and semi-analytical low-thrust \ac{CAM} methods}
    \label{tab:AnalyticalSummary_LowThrust}
    \begin{tabular}{llllll}
    \hline
        Ref.                     & Approach        & Constraint & Dynamics   \\ \hline
        \cite{Lee2014}                & Semi-Analytical & MD         & Keplerian  \\
        \cite{Reiter2018}             & Analytical      & MD         & Keplerian  \\
        \cite{Gonzalo2019Semi}        & Semi-Analytical & MD         & Drag + SRP \\
        \cite{Hernando-Ayuso2020}     & Analytical      & Chan       & Keplerian  \\
        \cite{DeVittori2022,DeVittori2022Geo,DeVittori2023}          & Analytical      & Chan, MD   & Keplerian  \\
        \cite{Strobel2023}             & Analytical      & Chan       & Keplerian  \\
        \cite{gonzalo2021computationally,gonzalo2022single,gonzalo37efficient} & Analytical & Chan & Keplerian \\
        \cite{DeVittori2023a}          & Semi-Analytical & Chan       & Keplerian  \\
        \cite{colombo2023sensitivity}  & Semi-Analytical & Chan       & Keplerian  \\
        \cite{DellElce2024,DeVeld2024} & Analytical      & MD         & Keplerian  \\

               \cite{Polli2025} & Semi-Analytical      & MD, SMD, Eclipse       & Keplerian  \\
       \hline
    \end{tabular}
\end{table}

In a consecutive effort to the impulsive formulation resumed in the previous section, Hernando-Ayuso and Bombardelli in \cite{Hernando-Ayuso2020} investigate optimal low-thrust \ac{CAM}s between satellites in circular orbits. In the study, the manoeuvre is divided into two consecutive parts: an initial firing arc, where the thrust magnitude is kept constant and the thrust direction is optimised, and a subsequent coasting phase, where the control is switched off until \ac{TCA}. After assuming negligible mass variations, Pontryagin’s Maximum Principle is leveraged to optimise the thrust direction while maintaining constant thrust magnitude. An indirect optimisation method, based on the hypotheses of short-term encounters and linearised dynamics, is presented in B-plane coordinates, along with a fully analytical solution for tangential manoeuvres. Research findings highlight that the orientation of the covariance ellipsoid can significantly impact the optimal timing of a low-thrust \ac{CAM}, with last-minute manoeuvrers being less sensitive, while also showing that non-Keplerian perturbations negligibly affect the problem's accuracy.

Building upon this, in \cite{DeVittori2022}, authors present analytical formulations for \ac{EO} and \ac{FO} \ac{CAM}s in short-term encounters, utilising constant, uncorrelated covariances and a spherical object approximation. The algorithms are developed in both \ac{ECI} and B-plane reference frames, with constraints enforced on \ac{MD}, \ac{SMD} and \ac{PoC}. Comparatively, \ac{EO} solutions in the B-plane exhibit slightly reduced accuracy due to the first-order approximation of dynamics. Optimal manoeuvrers, consistent with prior research, are typically near-tangential when planned in advance. Additionally, \ac{EO} solutions can be leveraged to solve \ac{FO} problems with bounded acceleration, achieving propellant efficiency. The algorithms were validated on realistic short-term conjunction scenarios in \ac{LEO}, demonstrating that the linearisations introduce negligible constraint violations, even for highly elliptical orbits, and that gravitational and atmospheric drag perturbations effects have a limited impact on the manoeuvres' safety and effectiveness.
A step forward is achieved in \cite{DeVittori2023} where previous algorithms are extended by incorporating return on the initial orbit with a point-to-point and a point-to-orbit strategy. In both policies, analytical solutions are retrieved by combining a series of \ac{EO} problems.
Finally, in the \ac{GEO} regime, \cite{DeVittori2022Geo} introduces algorithms that integrate low-thrust \ac{CAM}s and \ac{SK} by concurrently meeting a target \ac{PoC} and maximising the permanence time within the box. An \ac{EO} problem is devised as an analytical solution to a three-point boundary value problem, with conditions being the initial manoeuvring point, the target terminal state and the minimum \ac{PoC} at \ac{TCA}, and a linearisation performed through \ac{STM} linking initial and final states and costates. To ensure computational efficiency, Equinoctial Orbital Elements are employed and compared to a standard \ac{ECI} formulation. The solution can handle two cases: \ac{SK} alone, ensuring \ac{PoC} under a safety threshold; or, when not fulfilling this requirement, the algorithm detects the best strategy to command \ac{CAM} and \ac{SK} imposing an arbitrary \ac{PoC} limit at \ac{TCA} without redesigning the \ac{CAM}.

Gonzalo \textit{et al.} \cite{gonzalo2021computationally} introduce analytical and semi-analytical models for \ac{CAM}s using averaging techniques. This work, part of the MISS software tool for \ac{CAM} design, enhances the characterisation of phasing changes during close approaches, significantly improving accuracy with minimal added complexity. The model's efficiency is demonstrated through sensitivity analyses, and future improvements, such as the inclusion of normal thrust acceleration components, are discussed.
In a continuation of this work, an analytical framework for designing low-thrust CAMs under uncertainty using single-averaged dynamical models is presented in \cite{gonzalo2022single}. Building on the MISS software, authors incorporate both tangential and normal thrust components to optimise manoeuvre design. The method reduces \ac{PoC} using a quasi-optimal, piecewise constant thrust control profile derived from an eigenproblem approach. Moreover, this work models uncertainty using \ac{GMM}, which offers a practical solution for \ac{STM} in congested orbital environments. In their last contribution \cite{gonzalo37efficient}, the authors present analytical and semi-analytical models for low-thrust \ac{CAM} computation and design. With the aid of single-averaging of Keplerian motion equations, the secular and oscillatory components of orbit evolution due to a \ac{CAM} are effectively separated. For manoeuvre design, quasi-optimal piecewise-constant thrust profiles are derived based on their impulsive model counterparts.

Str{\"o}bel \textit{et al.} \cite{Strobel2023} propose an analytical method to determine the required coasting arc or burn duration to achieve a specified \ac{PoC}, modelled using Chan's series expansion. The method is based on the Hill-Clohessy-Wiltshire equations, which yield closed-form analytical solutions for the coasting arc following a short-duration manoeuvre. The approach assumes negligible perturbations, nearly circular orbits, minimal deviation between the object and the reference trajectory, and thrust applied in the along-track direction. The method effectively reduces \ac{PoC} when \ac{CAM}s are executed several revolutions before the \ac{TCA}, with execution time strongly influenced by encounter geometry. While achieving accuracy comparable to numerical simulations, the authors suggest its use as a preliminary estimate for more precise numerical methods.

Reiter \textit{et al.} \cite{Reiter2018} addresses the challenge of \ac{COLA} in \ac{LEO} under minimal notification times. It assumes scenarios with less than half an orbital period to plan and execute manoeuvrers, treating debris as a point mass with position uncertainty aligned along its velocity vector. The study recurs to finite burns with constant thrust magnitude to optimise thrust location, direction, and duration, aiming to minimise fuel usage while achieving safe separation, also expressed as \ac{PoC} threshold. A key contribution is the development of a semi-analytical formulation that links thrust duration, \ac{PoC}, and separation distance, facilitating rapid onboard decision-making. Case studies demonstrate that when notification times are below 20 minutes, reducing collision probability is prioritized over fuel consumption. Conversely, for longer notification times, minimizing manoeuvre duration becomes the primary objective.

Regarding semi-analytical techniques, Lee \textit{et al.} \cite{Lee2014} theorize how to generate near-optimal collision-free trajectories for active spacecraft. By incorporating a penalty function into the performance index, there is no need for iterative processes and initial guesses required in traditional optimisation. Using Pontryagin’s Principle and generating functions, the study derives feedback control laws as truncated power series, enabling efficient trajectory design.
Simulations demonstrate successful \ac{COLA} in single- and multi-obstacle scenarios with fuel efficiency close to reference solutions. While effective, it requires empirical tuning of penalty parameters and solving Hamilton-Jacobi equations. This approach offers a practical alternative to direct optimisation, especially in scenarios sensitive to computational limitations or initial conditions. Gonzalo \textit{et al.} \cite{Gonzalo2019Semi} develop a semi-analytical approach for \ac{CAM} analysis using proximal motion equations, incorporating effects such as drag and \ac{SRP}. Validated through a variety of test cases, the approach highlights its operational advantages.
The numerical efficiency of \ac{FO} \ac{CAM} policies is further enhanced in \cite{DeVittori2023a} with tangential thrust. The method addresses advanced conjunction notifications and balances computational efficiency with accuracy by leveraging an analytical propagator for dynamics. The algorithm starts by converting the acceleration profile derived from analytical \ac{EO} \ac{CAM} solution into a fuel-efficient bang-bang thrust profile. By iteratively refining the duration of firing windows, it achieves significant $\Delta v$ savings while ensuring compliance with operational constraints like thruster capabilities. Test cases demonstrate reduced computational times and practical applicability for both \ac{LEO} and \ac{GEO}, offering significant fuel savings with competitive performance \cite{colombo2023sensitivity}. 
Furthermore, in this work, the authors extend the approach to design \ac{EO} and \ac{FO} manoeuvres with radial thrust. In this case, \ac{DA} linearization serves for ascertaining the optimal switch-off time, given an ON-OFF profile and ensuring that the target \ac{SMD} is met. Ultimately, \cite{colombo2023sensitivity} introduces a novel approach for \ac{CAM} design during orbit-raising scenarios. The manoeuvre is achieved by determining the engine's optimal switch-off time, after which the subsequent drifting phase is analytically propagated until the anticipated \ac{TCA} is reached, with its duration being determined to ensure that the target \ac{SMD} is realised.

Dell'Elce, De Veld and Pomet in \cite{DellElce2024} and \cite{DeVeld2024} address the challenge of determining the minimum time required to increase \ac{MD} over a given threshold with a continuous propulsion up to \ac{TCA}. They propose an analytical approach using backward integration from \ac{TCA}, applying linearization principles and locally optimal steering laws. The first paper introduces a 7-dimensional state vector representation incorporating the time of closest approach and equinoctial elements: the necessary optimality conditions are derived using Pontryagin’s Maximum Principle at each time of the ballistic trajectory. An approximate solution is then proposed to estimate the minimum warning time to increase \ac{MD} over a predefined threshold; the method cannot handle \ac{PoC} constraints. The second paper focuses on deriving a locally optimal steering law validated through numerical simulations to assess its performance across various scenarios. In agreement with previous research, such as \cite{Hernando-Ayuso2020} and \cite{DeVittori2022}, the authors find that the thrust direction transitions from nearly purely radial near \ac{TCA} to predominantly tangential about half an orbit earlier. These studies offer valuable insights into the operational feasibility of low-thrust \ac{CAM}s under time-constrained scenarios, emphasizing the sensitivity of results to thrust magnitude and conjunction geometry. The computational efficiency of the method and its applicability to high-risk, low-warning time scenarios make it a promising candidate for future autonomous applications. 

Polli \textit{et al.} \cite{Polli2025} outline a semi-analytical model to incorporate eclipse constraints, preventing thrusting when in Earth's shadow, while simultaneously designing a fuel-optimal low-thrust \ac{CAM} that maximizes \ac{MD} or minimizes \ac{PoC}. The continuous low-thrust dynamics are discretized into a sequence of impulsive manoeuvres, allowing for a tractable representation of the otherwise continuous control problem. For each impulsive manoeuvre, the optimal thrust direction is determined via eigenvalue analysis of the quadratic form of the state transition matrix, which effectively captures the system's relative dynamics. Results from test cases show good agreement with high-fidelity simulations and the convergence reliability makes the method suitable for onboard implementation. The study also includes a sensitivity analysis, highlighting the model’s robustness under varying propulsion and lead-time scenarios.

%%%%%%%%%%%%%%%%%%%%%%%%%%%%%%%%%%%%%%%%%%%%%%%%%%%%%%%%%%%%%%%%%%%%%%%%%%%
% Comparison CAM Methods
%%%%%%%%%%%%%%%%%%%%%%%%%%%%%%%%%%%%%%%%%%%%%%%%%%%%%%%%%%%%%%%%%%%%%%%%%%%
\section{Comparison of Computationally Efficient CAM Routines in Single Short-term Encounters}
\label{sec:Comparison}

This section compares the six methods listed in \cref{tab:selectedmethods} that are best suited for autonomous applications.
To provide a comprehensive comparison, the algorithms are run using Keplerian dynamics on a set of test cases made accessible in the context of ESA's Collision Avoidance Challenge, available for download.\footnote{\url{github.com/arma1978/conjunction}} The set includes the \ac{CDM}s of 2.170 synthetic conjunctions in \ac{LEO}, with values of \ac{MD} spanning from $8.2$~\si{m} to $2$~\si{km} and \ac{PoC} between $10^{-6}$ and $0.136$. In \cref{tab:params}, the simulation parameters for \ac{CAM} design are showcased. Two decision metrics — \ac{MD} and \ac{PoC} — with four threshold levels are considered, given a warning time of one orbital period.

\begin{table}[tb!]
\centering
    \caption{Recap of the selected methods for the numerical comparison.}
    \label{tab:selectedmethods}
    \begin{tabular}{llcl}
        \hline
        Method & Ref.                      & Low-thrust       & Approach  \\  \hline
        Bomb15      & \cite{Bombardelli2015}         & \xmark            & Analytical  \\
        Arm21       & \cite{Armellin2021}            & \xmark            & SCP  \\
        DeVit23     & \cite{DeVittori2023a}          & \cmark           & Semi-Analytical  \\
        DeVit24     & \cite{DeVittori2024capitolo6}  & \xmark           & Analytical  \\
        Pav24a      & \cite{Pavanello2024Multiple}   & \cmark            & SCP  \\
        Pav24b      & \cite{Pavanello2024Recursive}  & \cmark            & PP  \\
        \hline
    \end{tabular}
\end{table}

\begin{table}[tb!]
    \caption{Comparison parameters.}
    \label{tab:params}
    \centering
    \begin{tabular}{ccc}
        \hline
        \ac{MD} thresholds & \ac{PoC} thresholds  & Warning time \\  \hline
        $[0.5, \ 1.5, \ 2.5, \ 3.5]$ & $[10^{-4}, \ 10^{-5}, \ 10^{-6}, \ 10^{-7}]$  & 1 revolution \\
        \hline
    \end{tabular}
\end{table}

The methods are evaluated in three of the four fundamental aspects for autonomy reported in \cref{sec:overallintroduction}: (i) efficiency, represented by the computation speed; (ii) optimality, which is connected to the $\dv$ required by the manoeuvre; (iii) reliability, i.e., the ability of to target the desired \ac{PoC} or \ac{MD} in a validation setup. 
Given that Bomb15 and DeVit24 need a specific firing time, a total of 10 possible manoeuvring points are considered, and the cheapest manoeuvre is selected. Pav24b instead is run using the \lq\lq nodes filtering" algorithm presented in the original paper, with 20 firing opportunities. The discretization for Pav24a and Arm21 consists of 60 firing opportunities. Lastly, DeVit23 solves the \ac{EO} problem before turning to the \ac{FO}, automatically selecting the best firing window. The maximum thrust available to the spacecraft is $10$ \si{mm/s^2}, consistent with a chemical system. The proposed simulation campaign does not allow for the evaluation of the robustness of the methods. The algorithms are run on 11th Gen Intel(R) Core(TM) i7-11700 @2.5GHz, 2496 Mhz, 8 Cores, 16 Logical Processors, using a combination of MATLAB R2024b, C++, Python, and Mosek 10.0.

\begin{figure}[tb!]
    \centering
    \subfloat[\ac{MD}.]{\includegraphics[width=0.48\textwidth]{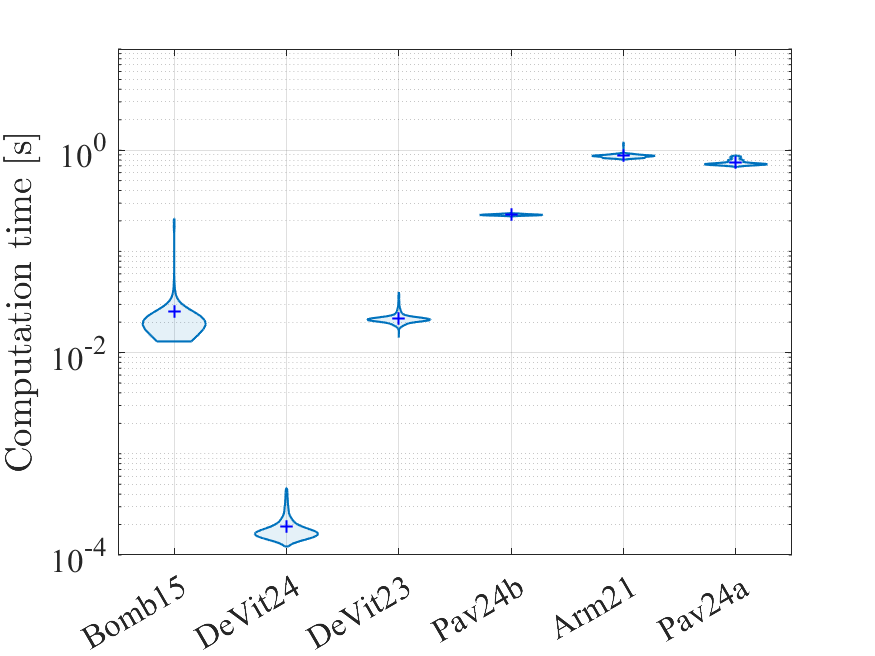}\label{fig:simTimeMd}}\hfill
    \subfloat[\ac{PoC}.]{\includegraphics[width=0.48\textwidth]{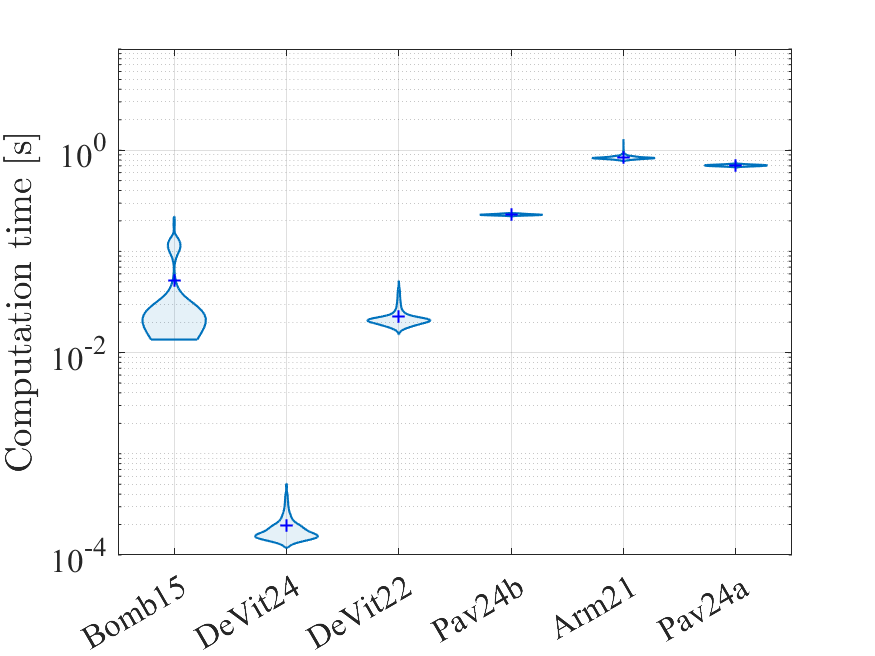}\label{fig:simTimePoC}}
    \caption{Distribution of simulation time.}
    \label{fig:simTime}
\end{figure}

In \cref{fig:simTime}, the distributions of computation times are shown in a logarithmic scale.
Using \ac{MD} or \ac{PoC} does not influence the runtimes significantly, with the only exception of Bomb15, which appears slightly slower in the \ac{PoC} case. The first three methods, being analytical or semi-analytical, are at least one order of magnitude faster than the numerical ones. In particular, the DeVit24 is the clear winner in this aspect, showing a distribution between $0.1$ \si{ms} and $1$ \si{ms}. 
The numerical approaches require higher computation times, mainly for the problem transcription and multiple runs of the interior-point method; however, in all simulations, the run time is kept below $1$ \si{s}, acceptable for autonomous use. The advantage of these methods, which is to include additional operational constraints, is not exploited in this setting.

\begin{figure}[tb!]
    \centering
    \subfloat[$\bar{d}_{miss} = 0.5$.]{\includegraphics[width=0.48\textwidth]{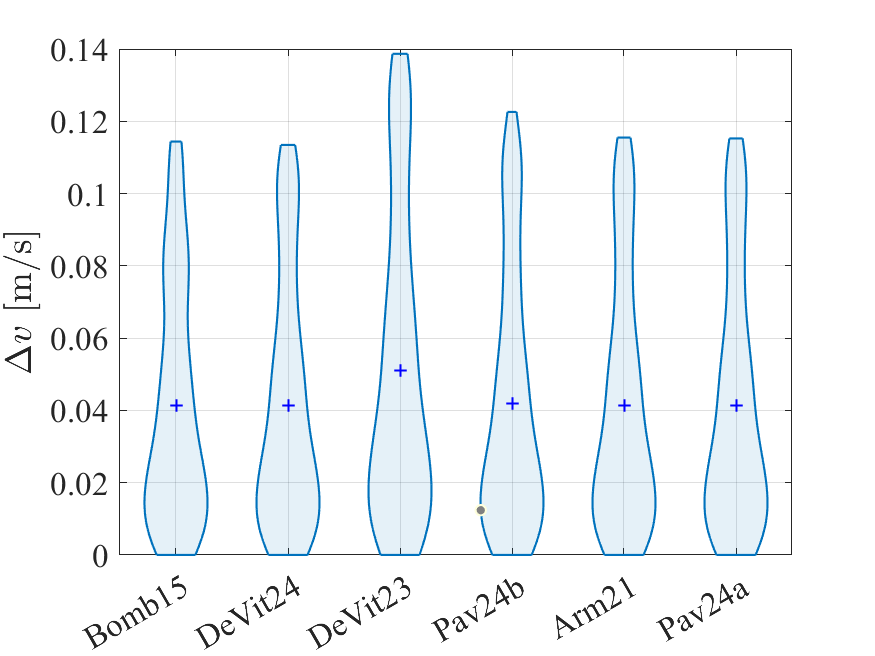}\label{fig:dvMd05}}\hfill
    \subfloat[$\bar{d}_{miss} = 1.5$.]{\includegraphics[width=0.48\textwidth]{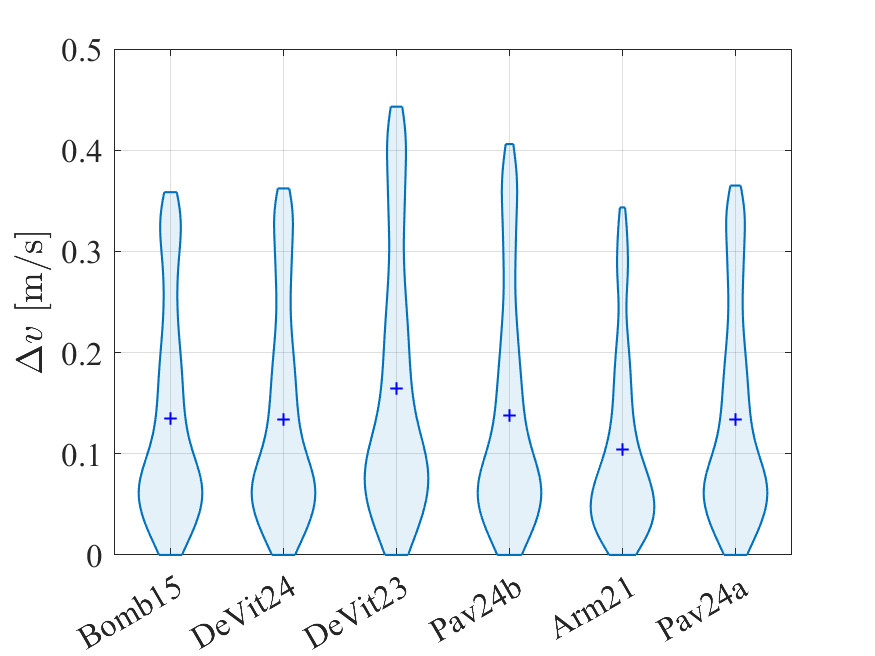}\label{fig:dvMd15}} \\
    \subfloat[$\bar{d}_{miss} = 2.5$.]{\includegraphics[width=0.48\textwidth]{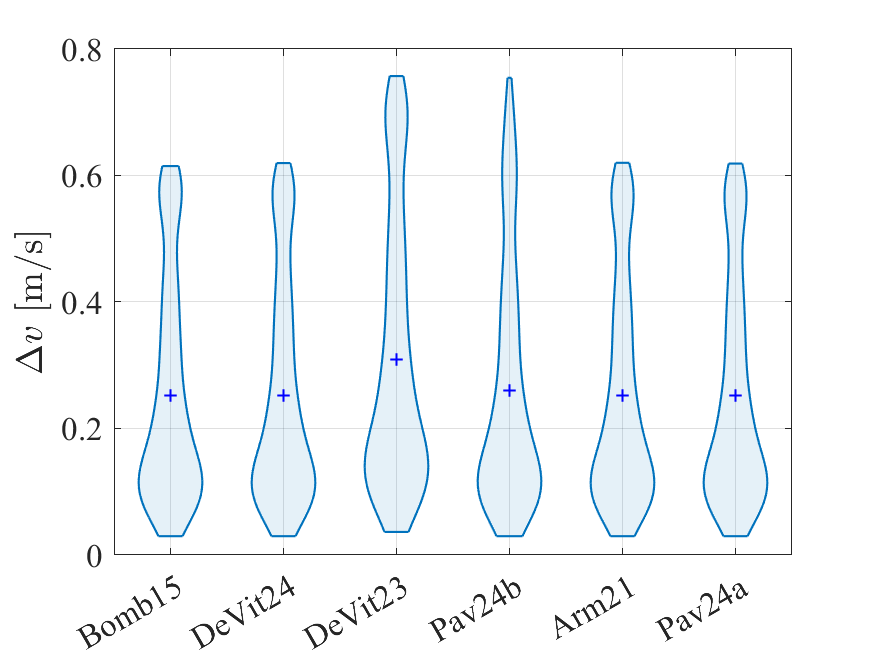}\label{fig:dvMd25}}\hfill
    \subfloat[$\bar{d}_{miss} = 3.5$.]{\includegraphics[width=0.48\textwidth]{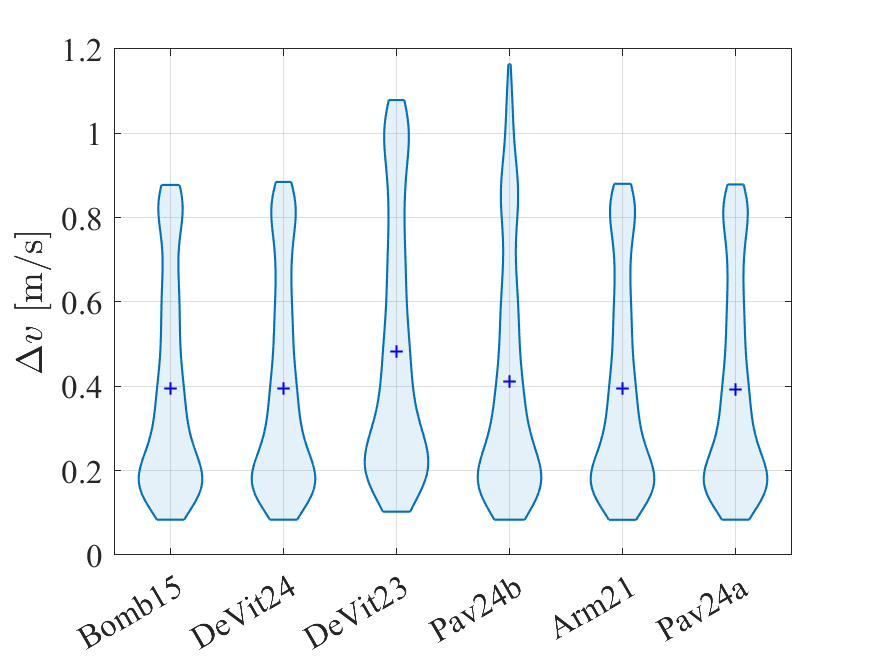}\label{fig:dvMd35}}
    \caption{\ac{MD}: Distribution of $\dv$s.}
    \label{fig:dvMD}
\end{figure}

\begin{figure}[tb!]
    \centering
    \subfloat[$\pclim=10^{-4}$.]{\includegraphics[width=0.48\textwidth]{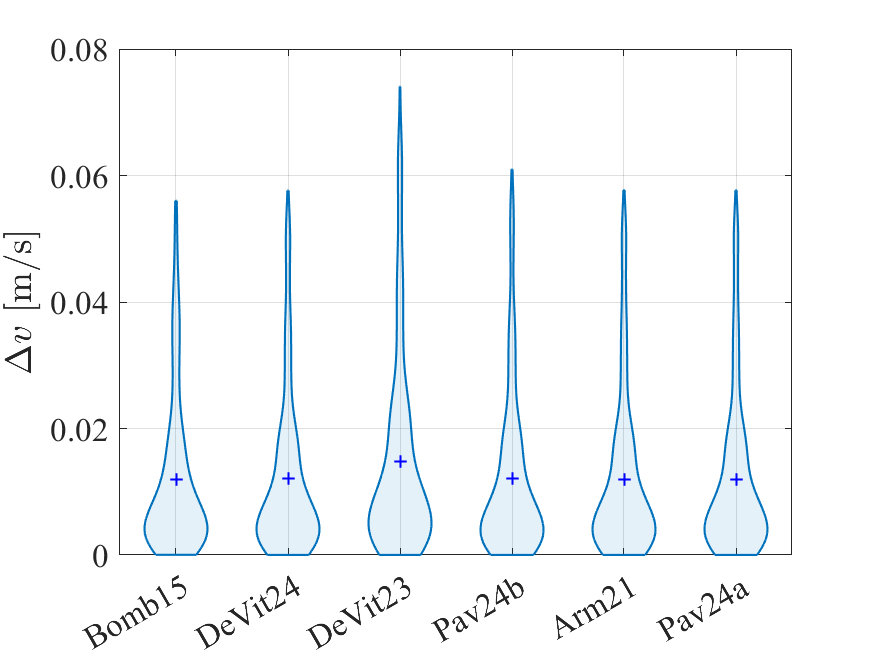}\label{fig:dvPoc4}}\hfill
    \subfloat[$\pclim=10^{-5}$.]{\includegraphics[width=0.48\textwidth]{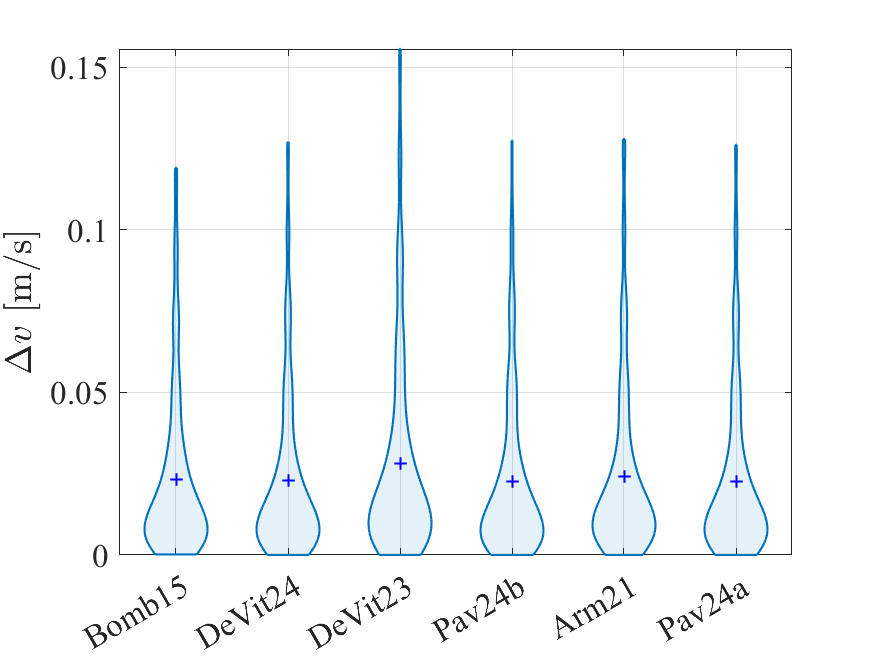}\label{fig:dvPoc5}} \\
    \subfloat[$\pclim=10^{-6}$.]{\includegraphics[width=0.48\textwidth]{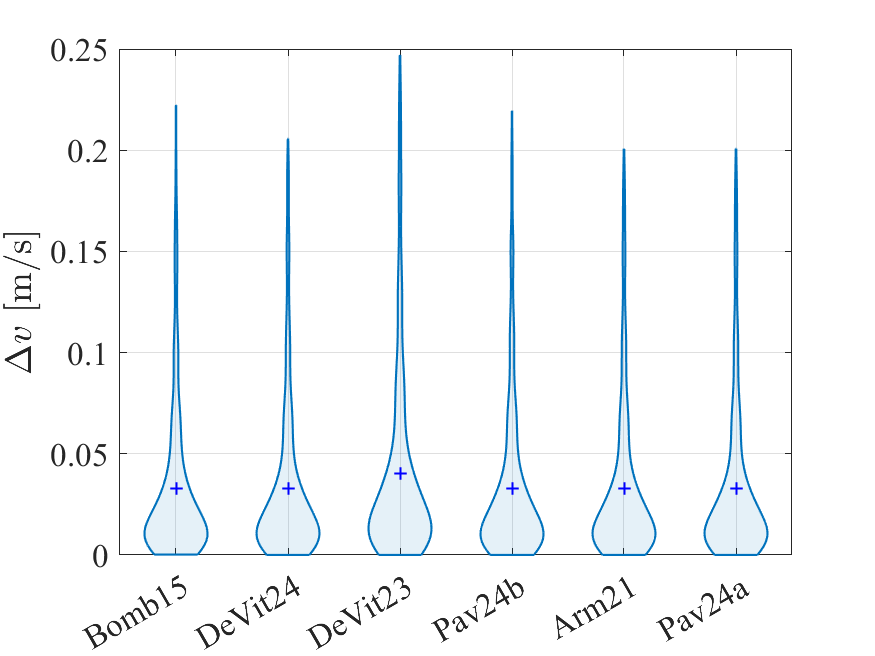}\label{fig:dvPoc6}}\hfill
    \subfloat[$\pclim=10^{-7}$.]{\includegraphics[width=0.48\textwidth]{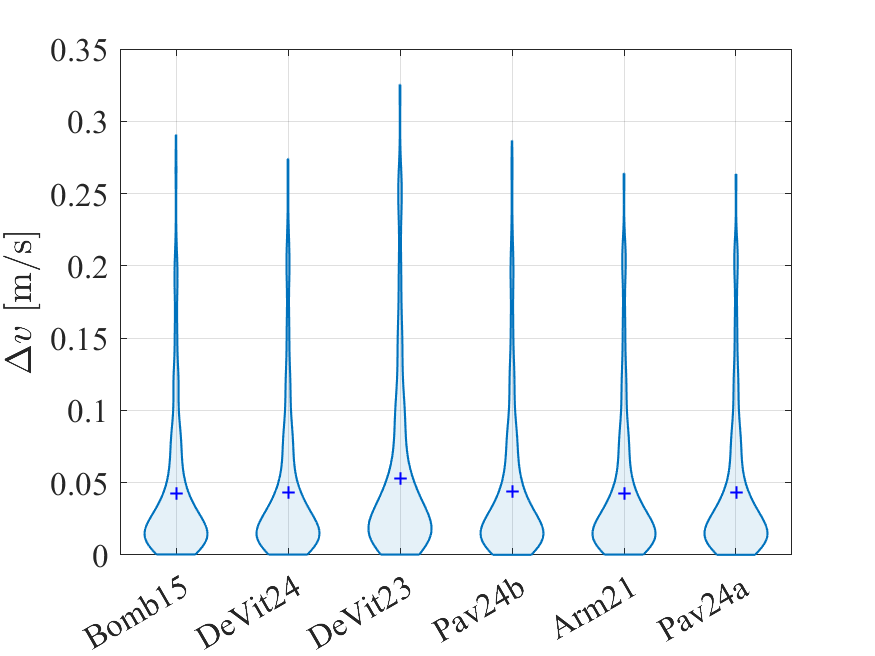}\label{fig:dvPoc7}}
    \caption{\ac{PoC}: Distribution of $\dv$s.}
    \label{fig:dvPoC}
\end{figure}

\begin{figure}[tb!]
    \centering
    \subfloat[\ac{MD}.]{\includegraphics[width=0.48\textwidth]{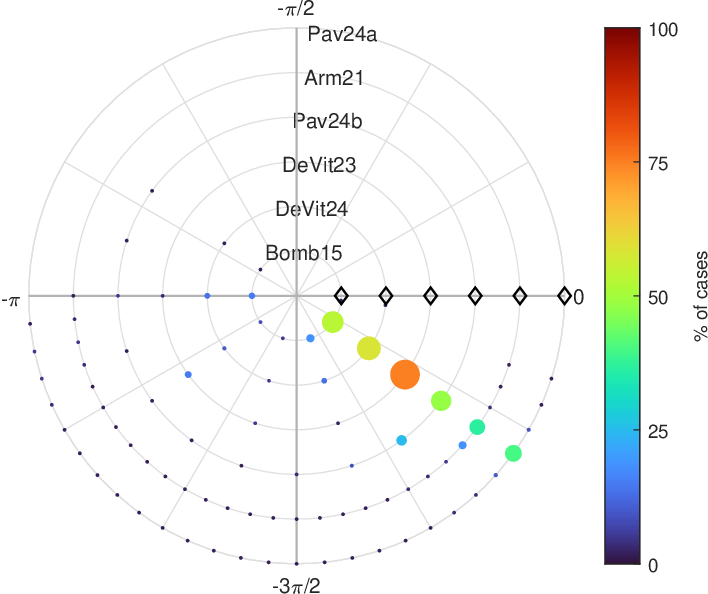}\label{fig:firingtimesmd}}\hfill
    \subfloat[\ac{PoC}.]{\includegraphics[width=0.48\textwidth]{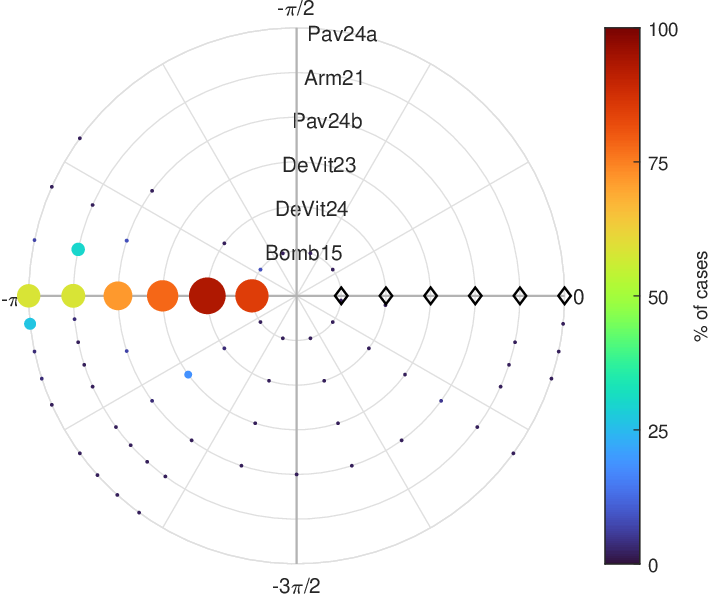}\label{fig:firingtimespoc}}
    \caption{Preferred firing times by each considered method. \ac{TCA} is represented by the rhomboid at 0; firing close to 0 counter-clockwise means that the manoeuvre happens close to one revolution before \ac{TCA}.}
    \label{fig:firingTimes}
\end{figure}

\cref{fig:dvMD} and \cref{fig:dvPoC} show the distributions of the computed manoeuvre magnitudes in terms of total $\Delta v$, while \cref{fig:firingTimes} depicts the manoeuvring points over a single orbit of warning time. Additionally, \cref{tab:meanDv} summarises the mean $\Delta v$ values for each threshold and method considered.
Comparing \cref{fig:dvMD} from (a) to (d), it is clear that the required $\Delta v$ for the \ac{CAM} is proportional to the \ac{MD} threshold. All methods exhibit a largely consistent distribution shape, suggesting that the manoeuvres they compute are similar. As expected, the low-trust approaches, like DeVit23, Pav24b, and Pav24a present a mean value (the blue cross) slightly higher than the impulsive methods. The mean value of DeVit23 is the highest in all cases because it selects either a purely tangential or a purely radial firing direction. This binary choice precludes the intermediate option — typically an almost fully tangential thrust with a minor radial component — that is available in the other implementations. Similar findings are observed for the $\Delta v$ distributions in the \ac{PoC} simulations in \cref{fig:dvPoC}. As before, DeVit23 indicates a slightly higher mean value, while the other techniques are closely aligned. Here, the $\Delta v$s are usually lower than those seen with the \ac{MD} threshold, and they gradually increase as the \ac{PoC} limit is lowered, as logically expected.

\begin{table}[tb!]
\caption{Mean value of the $\Delta V$ distribution for each considered metric threshold. The values are in \si{mm/s}.}
\label{tab:meanDv}
\centering
\begin{tabular}{l|cccccccc}
\hline
\multirow{2}{*}{Method} & \multicolumn{4}{c}{\ac{MD} limit [km]}              & \multicolumn{4}{|c}{\ac{PoC} limit [-]} \\ \cline{2-9} 
                        & 0.5   & 1.5   & 2.5   & 3.5 & \multicolumn{1}{|l}                    {1e-4}   & 1e-5    & 1e-6    & 1e-7    \\ \hline
Bomb15                  & 41.26 & 134.4 & 252.9 & \multicolumn{1}{l|}{395.1} & 11.98   & 23.38    & 32.73  & 42.55   \\
DeVit24                 & 41.34 & 133.9 & 252.4 & \multicolumn{1}{l|}{394.6} & 12.03   & 22.78    & 33.00  & 43.09   \\
DeVit23                 & 50.98 & 164.4 & 309.5 & \multicolumn{1}{l|}{483.2} & 14.80   & 28.03    & 40.50  & 52.79   \\
Pav24b                  & 41.89 & 137.8 & 260.7 & \multicolumn{1}{l|}{411.8} & 12.11   & 22.70    & 33.07  & 43.80   \\
Arm21                   & 41.37 & 133.9 & 251.8 & \multicolumn{1}{l|}{395.5} & 11.96   & 24.30    & 32.84  & 42.80   \\
Pav24a                  & 41.38 & 133.8 & 251.8 & \multicolumn{1}{l|}{393.6} & 12.00   & 22.75    & 32.91  & 42.90   \\ \hline
\end{tabular}
\end{table}

Many authors claim that the most efficient \ac{CAM}s are those performed when the spacecraft is in the opposite point on the orbit with respect to \ac{TCA}, i.e., $\pi$, $3\pi$, or $5\pi$, etc. true anomaly differences. This has been proven true for circular Keplerian orbits - so much so that operators adopt it as standard practice \cite{Zollo2024} - but it may not hold for even mildly elliptical orbits. This can be inferred by \cref{fig:firingtimespoc}, where the probability of a manoeuvre happening half a revolution before \ac{TCA} spans from  $60\%$ to $95$\% (depending on the method). One cannot neglect the possibility that a manoeuvre in a different portion of the orbit might be preferable, according to the results from the numerical methods. Nonetheless, as a general practice, it is safe to assume that manoeuvring half a revolution before \ac{TCA} is the best way to reduce \ac{PoC} in most applications.
However, as \cref{fig:firingtimesmd} shows, this is not the case when we want to increase \ac{MD}. All methods tend to agree that, as a general rule of thumb, it is more appropriate to execute the manoeuvre $324$ \si{deg} before \ac{TCA}. While the (semi-)analytical formulations assign over $50$\% of the manoeuvres to this slot, the methods with more refined grids - Pav24a and Arm21 - are in favour of a more widespread distribution over the course of the first half of the orbit. For these, only around $35$\% of the cases have an optimal manoeuvring point in $324$ \si{deg} before \ac{TCA}.

\begin{figure}[tb!]
    \centering
    \subfloat[\ac{MD}.]{\includegraphics[width=0.48\textwidth]{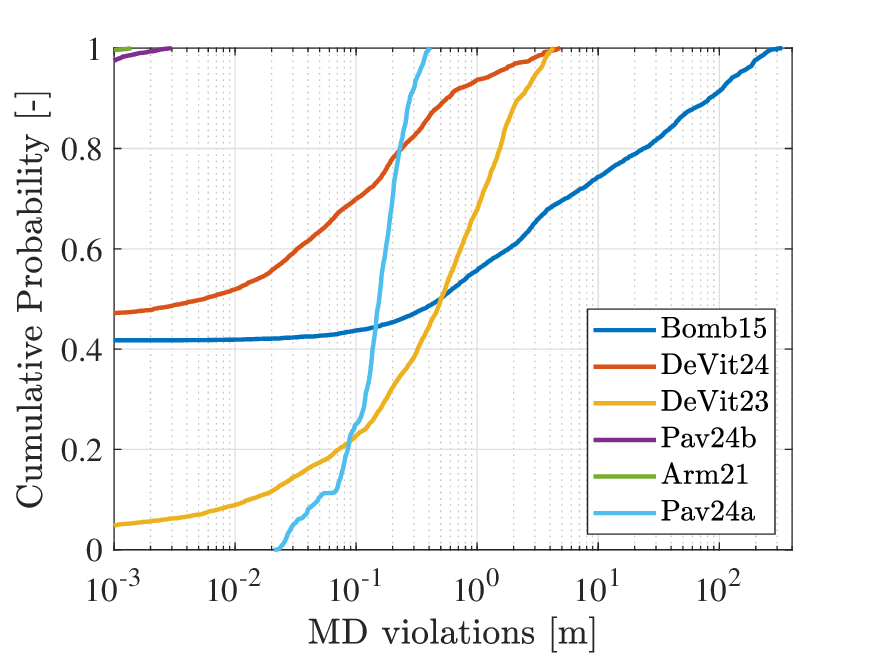}\label{fig:CDF_md}}\hfill
    \subfloat[\ac{PoC}.]{\includegraphics[width=0.48\textwidth]{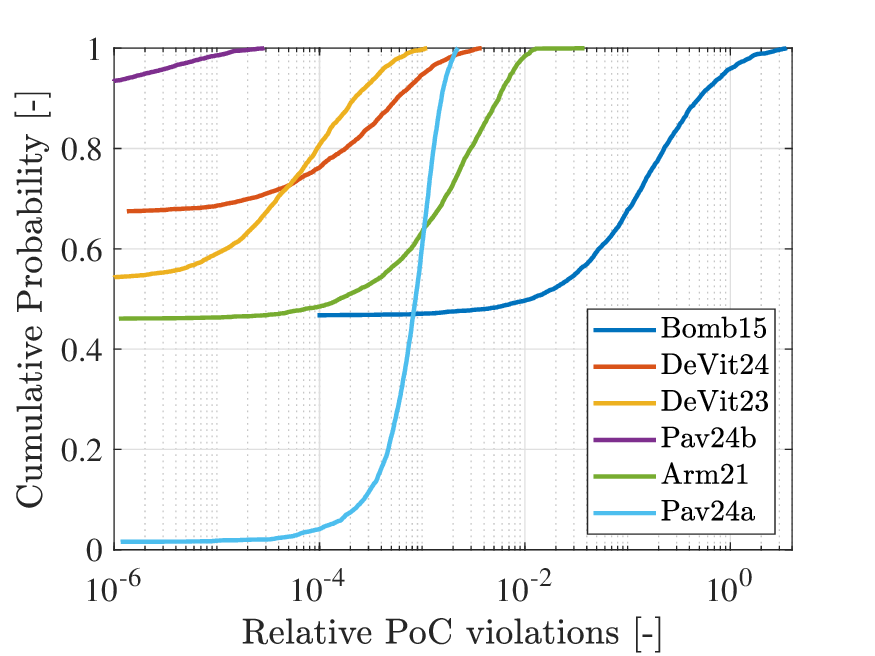}\label{fig:CDF_poc}}
    \caption{Constraint violations.}
    \label{fig:violations}
\end{figure}

Finally, reliability is evaluated considering the violations of the \ac{MD} or \ac{PoC} requirement. The impulses or the accelerations are applied at the specified times and fed to a Keplerian forward propagation scheme. The final position on the B-plane is used to validate \ac{MD} or \ac{PoC}. \ac{PoC} is computed with Chan's formula for all methods except Arm21, which employs Alfriend's one.
Given that \ac{MD} is a dimensional quantity, its violation is computed in terms of meters, whereas \ac{PoC} is evaluated in relative error; if \ac{PoC} is lower than the threshold -- or \ac{MD} is higher -- there is no violation, i.e., it is zero
\begin{equation*}
    \varepsilon_{P_C} = \max\left(0,\frac{\pc}{\pclim} - 1 \right), \quad \quad \quad
    \varepsilon_{\dmiss} = \max\left(0,\bar{d}_{miss} - \dmiss \right)
\end{equation*}
where $\varepsilon_{P_C}$ and $\varepsilon_{\dmiss}$ are the violations of \ac{PoC} and \ac{MD}, respectively, and $\pclim$ and $\bar{d}_{miss}$ are their corresponding thresholds.
In \cref{fig:violations}, the cumulative distribution functions of the violations are shown for the \ac{MD} simulations in \cref{fig:CDF_md} and for \ac{PoC} in \cref{fig:CDF_poc}. Bomb15 exceeds the \ac{MD} threshold by $100$ \si{m} in $10\%$ of cases and the \ac{PoC} one by more than $100\%$ in $5\%$ of cases. This may be due to Bomb15 underestimating the required $\Delta v$ compared to other methods, as shown in \cref{tab:meanDv}. All other methods keep errors below $1$ \si{m} and $3\%$ of the \si{PoC} threshold. Pav24b and Arm21 are the most accurate for \ac{MD}, with violations never exceeding a few \si{mm}. Pav24b tends to slightly overshoot, leading to a higher mean $\Delta v$ (see \cref{tab:meanDv}) but nearly negligible violations. Arm21 is less accurate in \ac{PoC} simulations, whereas Pav24a, DeVit23, and DeVit24 consistently remain within acceptable limits.

%%%%%%%%%%%%%%%%%%%%%%%%%%%%%%%%%%%%%%%%%%%%%%%%%%%%%%%%%%%%%%%%%%%%%%%%%%%
% Future Trends
%%%%%%%%%%%%%%%%%%%%%%%%%%%%%%%%%%%%%%%%%%%%%%%%%%%%%%%%%%%%%%%%%%%%%%%%%%%
\section{Future Trends}
\label{sec:FutureTrends}

This section explores recent \ac{CAM} advancements that we believe will become more and more significant in the future. In particular, \ac{CAM}s in the Cislunar domain and in thrustless applications (such as micro- and nanosats) are discussed.

\subsection{Cislunar Collision Avoidance}
\label{subsec:CislunarCA}

In response to growing concerns about space debris and orbital congestion, as highlighted in \cref{sec:overallintroduction}, both space agencies and private enterprises are increasingly investigating alternative orbital regimes for future space-based infrastructures. Among these, Cislunar space has garnered considerable attention thanks to its proximity, vast expanse, and critical scientific, strategic, and economic significance. As a result, numerous missions and infrastructures are being planned for this region. However, the challenges encountered in near-Earth orbits, such as the necessity for \ac{COLA}, may similarly emerge in this new orbital regime.

The literature on this topic is still in its early stages. In Ref. \cite{DeMaria2023}. the authors introduce an analytical method for impulsive \ac{CAM}s in Cislunar space using the \ac{CR3BP} model. A closed-form solution is then derived for the \ac{EO} formulation with free-firing direction, providing educated guess for a semi-analytical \ac{FO} \ac{CAM} with an ON-OFF profile. This last solution employs \ac{DA} for the optimisation of the firing window to retain some computational efficiency. Despite results typically obtained in two-body problems, resulting thrust profiles do not show any preference for the tangential direction. Ultimately, these algorithms adopt \ac{MD}, \ac{SMD} and \ac{PoC} as terminal constraints and are tested on a custom-designed test case in the L1-\ac{NRHO}, demonstrating highly promising computational times and accuracy. The same test case proposed in Ref. \cite{DeMaria2023} is used in Ref. \cite{Pavanello2024Recursive} to show that the recursive \ac{PP} implementation can deal with Cislunar dynamics, obtaining very similar results to the original work.

The challenges of conjunction assessment in the cislunar space is further developed by S. Crum et al. in \cite{Crum2024}, where authors evaluate collision risk and the need for \ac{CAM}s in this region. The results indicate that satellites in \ac{LLO} might require up to four manoeuvres annually to avoid collisions, with most close encounters occurring near the lunar equator. The analysis emphasises the significant impact of covariance matrix size on collision risk, stressing the need to improve the cislunar space domain awareness to mitigate future risks as lunar exploration increases.
Another work by T. Lechtenberg et al. \cite{lechtenberg2024novel} addresses the same topic by using internal cislunar population models and spatial-density approaches to estimate long-term collision risks. It also introduces techniques to ensure safety compliance despite gaps in state uncertainty data, tackling the increasing number of missions planned for the Moon and Cislunar space.

Ultimately, \cite{Nugent2025} leverage two dynamical models for \ac{CAM} design in this realm: \ac{CR3BP} and a \ac{HFEM}, both simulating spacecraft motion in the Earth-Moon system. The \ac{CR3BP} is used for initial analysis, while the \ac{HFEM} incorporates additional gravitational forces for greater accuracy. Collision avoidance is achieved through impulsive manoeuvres designed using the \ac{STM}, ensuring the spacecraft returns to the reference orbit after diversion, with evaluations conducted using both models.

%%%%%%%%%%%%%%%%%%%%%%%%%%%%%%%%%
\subsection{Thrustless CAM}
\label{subsec:AnalyticalCAM_Thrustless}

When performing a \ac{CAM}, the natural environmental forces present in orbit -- like \ac{SRP} or atmospheric drag -- are viable alternatives to actively thrusting. These solutions allow for changing the satellite's trajectory by increasing its ballistic coefficient or its \ac{SRP} equivalent area through attitude control or devices like drag/solar sails to maximise the environmental actions in the desired direction. The previously cited work in \cite{Gonzalo2019Semi} shows the possibility of employing \ac{SRP} and drag as a means to modify the orbit shape to avoid a close encounter, paving the way for \ac{CAM} capabilities also for small spacecraft by mounting solar or drag sails.

Mishne \textit{et al.} investigate in \cite{Mishne2017} the feasibility of drag and \ac{SRP} for \ac{COLA} in \ac{LEO} satellites. By orienting the satellite to maximise the effect of these natural forces, changes in the semi-major axis are achieved, creating sufficient along-track deviation to avoid collisions. It solely relies on attitude control and the difference between the satellite's maximum and minimum cross-sectional areas. Numerical simulations demonstrate that within three days, the technique can reduce \ac{PoC} below acceptable thresholds for satellites at altitudes between $600$ and $800$~\si{km}.
In \cite{Omar2020}, authors propose a way for \ac{COLA} based on controlling the spacecraft's ballistic coefficient to adjust its orbit. A simple analytical solution computes the difference in mean anomaly required to achieve a desired along-track miss distance, which is further refined numerically by modulating the ballistic coefficient using drag devices. The authors emphasise the importance of initiating manoeuvres early to minimise orbital deviations. Monte Carlo analysis demonstrates the robustness of the \ac{CAM} algorithm for operational scenarios, confirming the validity of the underlying assumptions. Building upon this work, F. Turco et al. in \cite{Turco2023} explore the use of atmospheric drag as an in-track \ac{COLA} strategy in \ac{LEO} by adjusting the satellite's ballistic coefficient through attitude changes. The method estimates the achievable separation distance between satellites experiencing varying drag forces, assuming circular orbits and a non-rotating atmosphere. The study incorporates constraints such as battery charging phases and conjunction geometry. The feasibility of the method is demonstrated using the Flying Laptop satellite, with results showing the feasibility of the method in reducing \ac{PoC}, despite some limitations connected to space weather and uncertainties.
Similarly, \cite{gaglio2024optimal} presents a novel \ac{COLA} strategy for satellites in \ac{VLEO}, leveraging aerodynamic drag modulation as a propellant-free alternative for unmanoeuvrable satellites. The approach balances two conflicting objectives -- increasing \ac{MD} while minimising orbital decay -- via a composed cost function. By employing the SOURCE CubeSat as a test case, the method is proved capable of achieving a $22$~\si{km} in-track separation with $165$~\si{m} orbital decay in a 24-hour warning period. Key contributions include a realistic model incorporating atmospheric dynamics and an \ac{OCP} ensuring practical control profiles. The study highlights its adaptability for future debris-congested environments and sustainable satellite operations.

Ultimately, in \cite{Kleinig2022} the study explores the use of ionospheric drag as a fuel-less mechanism for satellite \ac{COLA} in \ac{LEO}. By utilising an in-house propagation suite, the study investigates the impact of ionospheric and neutral drag on along-track separation for satellites at altitudes between $350$ and $500$~\si{km}. It further compares charged and uncharged methods, revealing that ionospheric drag can reduce the time to achieve NASA's approved collision risk reduction by at least 60\%. The results suggest that ionospheric drag offers a feasible alternative for \ac{COLA}, particularly for small satellites operating in \ac{LEO}.

%%%%%%%%%%%%%%%%%%%%%%%%%%%%%%%%%%%%%%%%%%%%%%%%%%%%%%%%%%%%%%%%%%%%%%%%%%%
% Conclusions
%%%%%%%%%%%%%%%%%%%%%%%%%%%%%%%%%%%%%%%%%%%%%%%%%%%%%%%%%%%%%%%%%%%%%%%%%%%
\section{Conclusions}
\label{sec:Conclusions}
This review provided a comprehensive analysis of computationally efficient collision avoidance manoeuvre (CAM) routines, focusing on their applicability in autonomous settings. Various numerical, semi-analytical, and analytical approaches have been examined, emphasizing their trade-offs in terms of computational speed, optimality, and reliability.
The study highlights that while heuristic and numerical optimization techniques offer flexibility in handling complex scenarios — including multiple conjunctions and operational constraints — they often require significant computational resources, making them less suitable for onboard autonomy. Conversely, analytical and semi-analytical methods provide rapid solutions but may lack the adaptability required for dynamically evolving conjunction assessment scenarios.

\section*{Acknowledgements}
The authors wish to thank Marco Felice Montaruli and Sergio Bonaccorsi for their precious contribution to the development of the codes used in the simulations presented in \cref{sec:Comparison}.

\bibliography{references}

\begin{thebibliography}{100}
\expandafter\ifx\csname url\endcsname\relax
  \def\url#1{\texttt{#1}}\fi
\expandafter\ifx\csname urlprefix\endcsname\relax\def\urlprefix{URL }\fi
\expandafter\ifx\csname href\endcsname\relax
  \def\href#1#2{#2} \def\path#1{#1}\fi

\bibitem{Kessler_syn}
D.~J. Kessler, B.~G. Cour-Palais, Collision frequency of artificial satellites: The creation of a debris belt, Journal of Geophysical Research: Space Physics 83~(A6) (1978) 2637--2646.
\newblock \href {http://dx.doi.org/https://doi.org/10.1029/JA083iA06p02637} {\path{doi:https://doi.org/10.1029/JA083iA06p02637}}.

\bibitem{ESA2024}
{ESA Space Debris Office}, {ESA's Annual Space Environment Report}, Document Type LOG GEN-DB-LOG-00288-OPS-SD, ESA ESOC, Robert-Bosch-Strasse 5, D-64293 Darmstadt, Germany, {ESA UNCLASSIFIED - Releasable to the Public. Status: Final.} (July 2024).

\bibitem{Horstmann2019}
A.~Horstmann, S.~Hesselbach, {ESA-MASTER: Hands-on}, in: {ESA-ECSL Space Debris Regulation, Standards, and Tools Workshop}, Vol.~1, ESA, 2019, p. n/a.

\bibitem{KosmosIridium}
N.~Johnson, The collision of iridium 33 and cosmos 2251: The shape of things to come, in: 59th International Astronautical Congress: IAC Proceedings, International Astronautical Federation (IAF), Seoul, Republic of Korea, 2009.

\bibitem{UNGUIDELINES}
{United Nations Office for Outer Space Affairs}, \href{https://www.un-ilibrary.org/content/books/9789213630921}{International Space Law}, United Nations, 2018.
\newline\urlprefix\url{https://www.un-ilibrary.org/content/books/9789213630921}

\bibitem{ESASpaceDebrisMitigationWorkingGroup2023}
{ESA Space Debris Mitigation Working Group}, {ESA Space Debris Mitigation Requirements}, Tech. rep., ESA, eSSB-ST-U-007 Issue 1 (2023).

\bibitem{DRAMA}
J.~Gelhaus, N.~Sanchez-Ortiz, V.~Braun, C.~Kebschull, J.~C. De~Oliveira, R.~Dominguez-Gonzalez, C.~Wiedemann, H.~Krag, P.~Vorsmann, Upgrade of drama-esa's space debris mitigation analysis tool suite, in: Proceedings of the 6th European Conference on Space Debris, Vol.~6, ESA, 2013.

\bibitem{AIAA2017}
M.~K. Jah, D.~Greiman, M.~Sengupta, S.~Magnus, P.~Melroy, S.~Helm, M.~Brown, {Space Traffic Management: Balancing Safety, Innovation, and Growth}, Tech. Rep. October, American Institute of Aeronautics and Astronautics (AIAA) (2017).

\bibitem{Krage2023}
F.~J. Krage, {Conjunction Assessment: Past, Present, and Future}, in: NASA Spacecraft Conjunction Assessment and Collision Avoidance Best Practices Handbook, NASA, Washington D.C., United States, 2023, Ch.~3, nASA/SP-20230002470 Rev 1.

\bibitem{CDM}
{Consultative Committee for Space Data Systems (CCSDS)}, \href{https://public.ccsds.org/Pubs/508x0b1e2s.pdf}{CCSDS Recommended Standard for Conjunction Data Messages}, issue 1 Edition, Vol. 508.0-B-1 of CCSDS Blue Book, CCSDS Secretariat, Space Communications and Navigation Office, NASA Headquarters, Washington, DC, USA, 2013.
\newline\urlprefix\url{https://public.ccsds.org/Pubs/508x0b1e2s.pdf}

\bibitem{CARA}
L.~K. Newman, A.~K. Mashiku, M.~D. Hejduk, M.~R. Johnson, J.~D. Rosa, Nasa conjunction assessment risk analysis (cara) updated requirements architecture, in: AAS/AIAA Astrodynamics Specialist Conference, American Institute of Aeronautics and Astronautics (AIAA), 2019.

\bibitem{aida2010collision}
S.~Aida, M.~Kirschner, M.~Wermuth, R.~Kiehling, Collision avoidance operations for leo satellites controlled by gsoc, in: Space Ops 2010, no. 2010-2298, German Space Operations Center DLR, Oberpfaffenhofen, Germany, 2010, presented at Space Ops 2010.

\bibitem{CONAN}
F.~Simarro, J.~Serrano, C.~Pérez, \href{https://www.eusst.eu/wp-content/uploads/2022/04/ICSSA-22_CONAN_v1.1.pdf}{Conan: Conjunction analysis tool}, Presentation at S3TOC-CONAN-ICSAA, Version 1.0 (March 2022).
\newline\urlprefix\url{https://www.eusst.eu/wp-content/uploads/2022/04/ICSSA-22_CONAN_v1.1.pdf}

\bibitem{SanchezOrtiz2014}
N.~S{\'{a}}nchez-Ortiz, K.~Merz, \href{https://www.researchgate.net/publication/261987324}{{CORAM: ESA's Collision Risk Assesment and Avoidance Manoeuvres Computation Tool}}, in: 2nd IAA Conference on Dynamics and Control of Space Systems, 2014.
\newline\urlprefix\url{https://www.researchgate.net/publication/261987324}

\bibitem{Cortesi2024}
S.~Patnala, A.~Abdin, {Spacecraft Collision Avoidance: Data Management, Risk Assessment, Decision Planning Models and Algorithms}, in: A.~Cortesi (Ed.), Space Data Management, Vol. 141 of Studies in Big Data, Springer Nature Singapore, Singapore, 2024, Ch.~2.
\newblock \href {http://dx.doi.org/10.1007/978-981-97-0041-7} {\path{doi:10.1007/978-981-97-0041-7}}.

\bibitem{Flohrer2022}
T.~Flohrer, K.~Merz, S.~Lemmens, B.~{Bastida Virgili}, J.~Siminski, F.~Letizia, Q.~Funke, F.~Mclean, {Current status and upcoming challenges in operational collision avoidance in a changing space debris environment}, in: 44th COSPAR Scientific Assembly, COSPAR, Athens, Greece, 2022.

\bibitem{Goldman2022}
D.~Goldman, \href{https://planet4589.org/astro/starsim/docs/Star2212.pdf}{{SpaceX Constellation Semi-Annual Report}}, Tech. rep., Space Exploration Technologies Corp., Washington, DC (2022).
\newline\urlprefix\url{https://planet4589.org/astro/starsim/docs/Star2212.pdf}

\bibitem{Gonzalo2020y}
J.~L. Gonzalo, C.~Colombo, P.~{Di Lizia}, {Introducing MISS, a new tool for collision avoidance analysis and design}, Journal of Space Safety Engineering 7~(3) (2020) 282--289.
\newblock \href {http://dx.doi.org/10.1016/j.jsse.2020.07.010} {\path{doi:10.1016/j.jsse.2020.07.010}}.

\bibitem{Gonzalo2021y}
J.~L. Gonzalo, C.~Colombo, {On-Board Collision Avoidance Applications Based on Machine Learning and Analytical Methods}, in: 8th European Conferences on Space Debris, no. May, ESA Space Debris Office, Darmstadt, Germany, 2021, pp. 20--23.

\bibitem{Sanchez2021}
L.~S{\'{a}}nchez, M.~Vasile, {Constrained Optimal Collision Avoidance Manoeuvre Allocation under Uncertainty for Subsequent Conjunction Events}, in: 72nd International Astronautical Congress: IAC Proceedings, International Astronautical Federation (IAF), Dubai, UAE, 2021.

\bibitem{Stroe2021}
I.~Stroe, A.~Stanculescu, P.~Ilioaica, C.~Blaj, M.~Nita, A.~Butu, D.~Escobar, J.~Tirado, B.~Bija, D.~Saez, {AUTOCA Autonomous Collision Avoidance System}, in: Proceedings of the 8th European Conference on Space Debris (virtual), no. May, ESA Space Debris Office, Darmstadt, Germany, 2021, p.~6.

\bibitem{Henry2023}
S.~Henry, R.~Armellin, T.~Gateau, {Safe-event pruning in spacecraft conjunction management}, Astrodynamics 7~(4) (2023) 401--413.
\newblock \href {http://dx.doi.org/10.1007/s42064-023-0165-5} {\path{doi:10.1007/s42064-023-0165-5}}.

\bibitem{Cannon2023Starling}
H.~Cannon, \href{https://ntrs.nasa.gov/citations/20230010760}{Starling swarm technology mission status and extended mission}, Small Satellite Conference NASA Short Talk, august 9, 2023, Document ID 20230010760 (2023).
\newline\urlprefix\url{https://ntrs.nasa.gov/citations/20230010760}

\bibitem{BastidaVirgili2019}
B.~{Bastida Virgili}, T.~Flohrer, H.~Krag, K.~Merz, S.~Lemmens, {CREAM - ESA's Proposal for Collision Risk Estimation and Automated Mitigation}, in: First International Orbital Debris Conference, Sugar Land, TX, 2019.

\bibitem{maybeck1982_book}
P.~Maybeck, \href{https://books.google.it/books?id=L_YVMUJKNQUC}{Stochastic Models, Estimation, and Control}, Elsevier Science, 1982.
\newline\urlprefix\url{https://books.google.it/books?id=L_YVMUJKNQUC}

\bibitem{Fuller_art}
A.~Fuller, Analysis of nonlinear stochastic systems by means of the fokker–planck equation, International Journal of Control 9 (2007) 603--655.
\newblock \href {http://dx.doi.org/10.1080/00207176908905786} {\path{doi:10.1080/00207176908905786}}.

\bibitem{Kumar}
M.~Kumar, S.~Chakravorty, J.~Junkins, A semianalytic meshless approach to the transient fokker–planck equation, Probabilistic Engineering Mechanics 25 (2010) 323--331.
\newblock \href {http://dx.doi.org/10.1016/j.probengmech.2010.01.006} {\path{doi:10.1016/j.probengmech.2010.01.006}}.

\bibitem{Sun_prop}
Y.~Sun, M.~Kumar, Uncertainty propagation in orbital mechanics via tensor decomposition, Celestial Mechanics and Dynamical Astronomy 124~(3) (2016) 269–--294.
\newblock \href {http://dx.doi.org/10.1007/s10569-015-9662-z} {\path{doi:10.1007/s10569-015-9662-z}}.

\bibitem{Junkins1996NonGaussianEP_article}
J.~L. Junkins, M.~R. Akella, K.~T. Alfriend, \href{https://api.semanticscholar.org/CorpusID:125478821}{Non-gaussian error propagation in orbital mechanics}, Journal of The Astronautical Sciences 44 (1996) 541--563.
\newline\urlprefix\url{https://api.semanticscholar.org/CorpusID:125478821}

\bibitem{dagumMC}
P.~Dagum, R.~Karp, M.~Luby, S.~Ross, An optimal algorithm for monte carlo estimation, in: Proceedings of IEEE 36th Annual Foundations of Computer Science, 1995, pp. 142--149.
\newblock \href {http://dx.doi.org/10.1109/SFCS.1995.492471} {\path{doi:10.1109/SFCS.1995.492471}}.

\bibitem{morselli2015high}
A.~Morselli, R.~Armellin, P.~{Di Lizia}, F.~{Bernelli Zazzera}, A high order method for orbital conjunctions analysis: Monte carlo collision probability computation, Advances in Space Research 55~(1) (2015) 311--333.
\newblock \href {http://dx.doi.org/https://doi.org/10.1016/j.asr.2014.09.003} {\path{doi:https://doi.org/10.1016/j.asr.2014.09.003}}.

\bibitem{Battin1987AnIT}
R.~H. Battin, An introduction to the mathematics and methods of astrodynamics, American Institute of Aeronautics and Astronautics (AIAA), 1987.
\newblock \href {http://dx.doi.org/10.2514/4.861543} {\path{doi:10.2514/4.861543}}.

\bibitem{RAO}
P.~P. Rao, S.~C. Bell, Conditional performance error covariance analyses for commercial titan launch vehicles, Journal of Guidance, Control, and Dynamics 14~(2) (1991) 398--405.
\newblock \href {http://dx.doi.org/10.2514/3.20652} {\path{doi:10.2514/3.20652}}.

\bibitem{gelb1974}
A.~Gelb, Applied Optimal Estimation, The MIT Press, 1974.

\bibitem{Julier}
S.~Julier, J.~Uhlmann, Reduced sigma point filters for the propagation of means and covariances through nonlinear transformations, in: Proceedings of the 2002 American Control Conference, Vol.~2, 2002, pp. 887--892 vol.2.
\newblock \href {http://dx.doi.org/10.1109/ACC.2002.1023128} {\path{doi:10.1109/ACC.2002.1023128}}.

\bibitem{Arasaratnam}
I.~Arasaratnam, S.~Haykin, Cubature kalman filters, IEEE Transactions on Automatic Control 54~(6) (2009) 1254--1269.
\newblock \href {http://dx.doi.org/10.1109/TAC.2009.2019800} {\path{doi:10.1109/TAC.2009.2019800}}.

\bibitem{Park}
R.~Park, D.~Scheeres, Nonlinear mapping of gaussian statistics: Theory and applications to spacecraft trajectory design, Journal of Guidance Control and Dynamics 29~(11) (2006) 1367--1375.
\newblock \href {http://dx.doi.org/10.2514/1.20177} {\path{doi:10.2514/1.20177}}.

\bibitem{Berz1999ModernMM_PDF}
B.~Martin, Chapter 1 - dynamics of particles and fields, in: P.~Hawkes (Ed.), Modern Map Methods in Particle Beam Physics, Vol. 108 of Advances in Imaging and Electron Physics, Elsevier, 1999, pp. 1--79.
\newblock \href {http://dx.doi.org/doi.org/10.1016/S1076-5670(08)70227-1} {\path{doi:doi.org/10.1016/S1076-5670(08)70227-1}}.

\bibitem{Massari2017}
M.~Massari, P.~{Di Lizia}, M.~Rasotto, {Nonlinear Uncertainty Propagation in Astrodynamics Using Differential Algebra and Graphics Processing Units}, Journal of Aerospace Information Systems 14~(9) (2017) 493--503.
\newblock \href {http://dx.doi.org/10.2514/1.I010535} {\path{doi:10.2514/1.I010535}}.

\bibitem{Sun2019}
Z.-J. Sun, Y.-Z. Luo, P.~di~Lizia, F.~B. Zazzera, {Nonlinear orbital uncertainty propagation with differential algebra and Gaussian mixture model}, Science China Physics, Mechanics \& Astronomy 62~(3) (2019) 34511.
\newblock \href {http://dx.doi.org/10.1007/s11433-018-9267-6} {\path{doi:10.1007/s11433-018-9267-6}}.

\bibitem{Crisan}
D.~Crisan, J.~M{\'i}guez, {Particle-kernel estimation of the filter density in state-space models}, Bernoulli 20~(4) (2014) 1879 -- 1929.
\newblock \href {http://dx.doi.org/10.3150/13-BEJ545} {\path{doi:10.3150/13-BEJ545}}.

\bibitem{Vittaldev2016Space}
V.~Vittaldev, R.~P. Russell, {Space object collision probability using multidirectional Gaussian mixture models}, Journal of Guidance, Control, and Dynamics 39~(9) (2016) 2161--2167.
\newblock \href {http://dx.doi.org/10.2514/1.G001610} {\path{doi:10.2514/1.G001610}}.

\bibitem{maestrini2023electrocam}
M.~Maestrini, A.~De~Vittori, J.~L. Gonzalo~Gómez, C.~Colombo, P.~Di~Lizia, J.~Míguez~Arenas, M.~Sanjurjo~Rivo, A.~Díez~Martín, P.~Gago~Padreny, D.~Escobar~Antón, Electrocam: assessing the effect of low-thrust uncertainties on orbit propagation, in: 2nd NEO and Debris Detection Conference, ESA Space Debris Office, ESA/ESOC, Darmstadt, Germany, 2023.

\bibitem{Alfano2022_SatOperatorsSurvey}
S.~Alfano, D.~L. Oltrogge, L.~Arona, Operators’ requirements for ssa services, Journal of the Astronautical Sciences 69 (2022) 1441--1476.
\newblock \href {http://dx.doi.org/10.1007/s40295-022-00346-8} {\path{doi:10.1007/s40295-022-00346-8}}.

\bibitem{Chan2008Spacecraft}
K.~Chan, {Spacecraft Collision Probability}, The Aerospace Press, El Segundo, USA, 2008.
\newblock \href {http://dx.doi.org/10.1017/S0001924000087546} {\path{doi:10.1017/S0001924000087546}}.

\bibitem{Chan_short}
K.~Chan, Comparison of methods for spacecraft collision probability computations, in: International Space Conference of Pacific-basin Societies (ISCOPS), Montreal, Canada, 2020.

\bibitem{Armellin2021}
R.~Armellin, {Collision avoidance maneuver optimization with a multiple-impulse convex formulation}, Acta Astronautica 186 (2021) 347--362.
\newblock \href {http://arxiv.org/abs/2101.07403} {\path{arXiv:2101.07403}}, \href {http://dx.doi.org/10.1016/j.actaastro.2021.05.046} {\path{doi:10.1016/j.actaastro.2021.05.046}}.

\bibitem{Foster1992}
J.~L. Foster, H.~S. Estes, {A Parametric Analysis of Orbital Debris Collision Probability and Maneuver Rate for Space Vehicles}, Tech. rep. (1992).

\bibitem{Alfriend2000}
K.~T. Alfriend, M.~R. Akella, J.~Frisbee, J.~L. Foster, D.-J. Lee, M.~Wilkins, {Probability of collision error analysis}, Space Debris 1 (1999) 21--35.
\newblock \href {http://dx.doi.org/https://doi.org/10.1023/A:1010056509803} {\path{doi:https://doi.org/10.1023/A:1010056509803}}.

\bibitem{Chan2004Short}
K.~Chan, Short-term vs. long-term spacecraft encounters, in: AAS/AIAA Astrodynamics Specialist Conference and Exhibit, American Institute of Aeronautics and Astronautics (AIAA).
\newblock \href {http://arxiv.org/abs/https://arc.aiaa.org/doi/pdf/10.2514/6.2004-5460} {\path{arXiv:https://arc.aiaa.org/doi/pdf/10.2514/6.2004-5460}}, \href {http://dx.doi.org/10.2514/6.2004-5460} {\path{doi:10.2514/6.2004-5460}}.

\bibitem{Chan2004International}
K.~Chan, International space station collision probability, in: AIAA/AAS Astrodynamics Specialist Conference, 2008, pp. 307--314.
\newblock \href {http://dx.doi.org/10.2514/6.2008-6774} {\path{doi:10.2514/6.2008-6774}}.

\bibitem{Bai2010Space}
X.-Z. Bai, L.~Chen, G.-J. Tang, {Space objects maximum collision probability analysis based on explicit expression}, Yuhang Xuebao/Journal of Astronautics 31~(3) (2010) 880 -- 887.
\newblock \href {http://dx.doi.org/10.3873/j.issn.1000-1328.2010.03.042} {\path{doi:10.3873/j.issn.1000-1328.2010.03.042}}.

\bibitem{Serra2016Fast}
R.~Serra, D.~Arzelier, M.~Joldes, J.-B. Lasserre, A.~Rondepierre, B.~Salvy, {Fast and accurate computation of orbital collision probability for short-term encounters}, Journal of Guidance, Control, and Dynamics 39~(5) (2016) 1009 -- 1021.
\newblock \href {http://dx.doi.org/10.2514/1.G001353} {\path{doi:10.2514/1.G001353}}.

\bibitem{Serra2015}
R.~Serra, D.~Arzelier, M.~Joldes, A.~Rondepierre, {Probabilistic Collision Avoidance for Long-term Space Encounters via Risk Selection}, in: Advances in Aerospace Guidance, Navigation and Control, Springer International Publishing, 2015, pp. 679--698.
\newblock \href {http://dx.doi.org/10.1007/978-3-319-17518-839} {\path{doi:10.1007/978-3-319-17518-839}}.

\bibitem{Patera2001General}
R.~P. Patera, {General method for calculating satellite collision probability}, Journal of Guidance, Control, and Dynamics 24~(4) (2001) 716--722.
\newblock \href {http://dx.doi.org/10.2514/2.4771} {\path{doi:10.2514/2.4771}}.

\bibitem{Alfano2005Numerical}
S.~Alfano, {A Numerical Implementation of Spherical Object Coliision Probability}, The Journal of the Astronautical Sciences 53~(1) (2005) 103--109.

\bibitem{Pavanello2024}
Z.~Pavanello, L.~Pirovano, R.~Armellin, {Long-Term Fuel-Optimal Collision Avoidance Maneuvers with Station-Keeping Constraints}, Journal of Guidance, Control, and Dynamics 47~(9) (2024) 1855--1871.
\newblock \href {http://dx.doi.org/10.2514/1.G007839} {\path{doi:10.2514/1.G007839}}.

\bibitem{Li2022}
J.~S. Li, Z.~Yang, Y.~Z. Luo, {A review of space-object collision probability computation methods}, Astrodynamics 6~(2) (2022) 95--120.
\newblock \href {http://dx.doi.org/10.1007/s42064-021-0125-x} {\path{doi:10.1007/s42064-021-0125-x}}.

\bibitem{Zhang2020}
S.~Zhang, T.~Fu, D.~Chen, H.~Cao, {Satellite instantaneous collision probability computation using equivalent volume cuboids}, Journal of Guidance, Control, and Dynamics 43~(9) (2020) 1757--1763.
\newblock \href {http://dx.doi.org/10.2514/1.G004711} {\path{doi:10.2514/1.G004711}}.

\bibitem{Patera2003}
R.~P. Patera, G.~E. Peterson, {Space vehicle maneuver method to lower collision risk to an acceptable level}, Journal of Guidance, Control, and Dynamics 26~(2) (2003) 233--237.
\newblock \href {http://dx.doi.org/10.2514/2.5063} {\path{doi:10.2514/2.5063}}.

\bibitem{Pineiro2011}
J.~J. Pi{\~{n}}eiro, E.~M. Romero, J.~Fuentes, {Optimal collision avoidance maneuvers using Pseudospectral Methods}, in: European Space Surveillance Conference, European Space Agency (ESA), 2011.

\bibitem{Lee2012Collision}
S.~C. Lee, H.~D. Kim, J.~Suk, {Collision avoidance maneuver planning using GA for LEO and GEO satellite maintained in keeping area}, International Journal of Aeronautical and Space Sciences 13~(4) (2012) 474--483.
\newblock \href {http://dx.doi.org/10.5139/IJASS.2012.13.4.474} {\path{doi:10.5139/IJASS.2012.13.4.474}}.

\bibitem{Greco2021}
C.~Greco, L.~S{\'{a}}nchez, M.~Manzi, M.~Vasile, {A Robust Bayesian Agent for Optimal Collision Avoidance Manoeuvre Planning}, in: Proceedings of the 8th European Conference on Space Debris - ESA/ESOC, no. April, Darmstadt, Germany, 2021, pp. 1--11.

\bibitem{Dutta2022}
S.~Dutta, A.~K. Misra, {Convex optimization of collision avoidance maneuvers in the presence of uncertainty}, Acta Astronautica 197 (2022) 257--268.
\newblock \href {http://dx.doi.org/10.1016/j.actaastro.2022.05.038} {\path{doi:10.1016/j.actaastro.2022.05.038}}.

\bibitem{Dutta2023a}
S.~Dutta, A.~K. Misra, {Effect of the nature of uncertainty on the optimization of collision avoidance maneuvers}, Advances in Space Research 72~(10) (2023) 4132--4146.
\newblock \href {http://dx.doi.org/10.1016/j.asr.2023.08.012} {\path{doi:10.1016/j.asr.2023.08.012}}.

\bibitem{SANCHEZ20232627}
L.~Sánchez, M.~Vasile, Intelligent decision support for collision avoidance manoeuvre planning under uncertainty, Advances in Space Research 72~(7) (2023) 2627--2648.
\newblock \href {http://dx.doi.org/doi.org/10.1016/j.asr.2022.09.023} {\path{doi:doi.org/10.1016/j.asr.2022.09.023}}.

\bibitem{Pavanello2024Multiple}
Z.~Pavanello, L.~Pirovano, R.~Armellin, A.~{De Vittori}, P.~{Di Lizia}, {A Convex Optimization Method for Multiple Encounters Collision Avoidance Maneuvers}, in: AIAA SciTech Forum, no. January, American Institute of Aeronautics and Astronautics (AIAA), Orlando, Florida, 2024, pp. 1--17.
\newblock \href {http://dx.doi.org/10.2514/6.2024-0845} {\path{doi:10.2514/6.2024-0845}}.

\bibitem{Pavanello2024LowThrust}
Z.~Pavanello, L.~Pirovano, R.~Armellin, A.~{De Vittori}, P.~{Di Lizia}, {A Convex Formulation for Collision Avoidance Maneuver Strategies During Low-Thrust Phases}, in: AIAA SciTech Forum, no. January, American Institute of Aeronautics and Astronautics (AIAA), Orlando, Florida, 2024, pp. 1--19.
\newblock \href {http://dx.doi.org/10.2514/6.2024-0844} {\path{doi:10.2514/6.2024-0844}}.

\bibitem{Pavanello2024Recursive}
Z.~Pavanello, L.~Pirovano, R.~Armellin, Recursive polynomial method for fast collision avoidance maneuver design, IEEE Transactions on Aerospace and Electronic Systems (2024) 1--12\href {http://dx.doi.org/10.1109/TAES.2024.3480041} {\path{doi:10.1109/TAES.2024.3480041}}.

\bibitem{Lopez2024}
D.~Pérez~López, A.~De~Vittori, P.~Di~Lizia, D.~Mortari, Exploring functional connections theory and linearized approaches in collision avoidance maneuver design: A comparative study, in: 75th International Astronautical Conference (IAC), 2024, pp. 1--16.
\newblock \href {http://dx.doi.org/10.52202/078360-0162} {\path{doi:10.52202/078360-0162}}.

\bibitem{Pavanello2024IAC}
Z.~Pavanello, L.~Pirovano, R.~Armellin, {Polynomial Recursive Method for Fast Collision Avoidance Manoeuvre with Station-Keeping}, in: 75th International Astronautical Congress (IAC), Milan, Italy, 2024, pp. 1--15.

\bibitem{Liu2017}
X.~Liu, P.~Lu, B.~Pan, {Survey of convex optimization for aerospace applications}, Astrodynamics 1~(1) (2017) 23--40.
\newblock \href {http://dx.doi.org/10.1007/s42064-017-0003-8} {\path{doi:10.1007/s42064-017-0003-8}}.

\bibitem{Malyuta2021Tutorial}
D.~Malyuta, T.~P. Reynolds, M.~Szmuk, T.~Lew, R.~Bonalli, M.~Pavone, B.~A{\c{c}}ıkmeşe, {Convex Optimization for Trajectory Generation: A Tutorial on Generating Dynamically Feasible Trajectories Reliably and Efficiently}, IEEE Control Systems 42~(5) (2022) 40--113.
\newblock \href {http://arxiv.org/abs/2106.09125} {\path{arXiv:2106.09125}}, \href {http://dx.doi.org/10.1109/MCS.2022.3187542} {\path{doi:10.1109/MCS.2022.3187542}}.

\bibitem{Wang2024}
Z.~Wang, {A survey on convex optimization for guidance and control of vehicular systems}, Annual Reviews in Control 57~(April) (2024) 100957.
\newblock \href {http://dx.doi.org/10.1016/j.arcontrol.2024.100957} {\path{doi:10.1016/j.arcontrol.2024.100957}}.

\bibitem{Boyd2}
S.~Boyd, L.~Vandenberghe, {Convex Optimization}, Cambridge University Press, Cambridge, NY, 2004.

\bibitem{Mao2017}
Y.~Mao, D.~Dueri, M.~Szmuk, B.~Açıkmeşe, Successive convexification of non-convex optimal control problems with state constraints, IFAC-PapersOnLine 50~(1) (2017) 4063--4069, 20th IFAC World Congress.
\newblock \href {http://dx.doi.org/doi.org/10.1016/j.ifacol.2017.08.789} {\path{doi:doi.org/10.1016/j.ifacol.2017.08.789}}.

\bibitem{Mao2021}
Y.~Mao, B.~Acikmese, {SCvx-fast: A Superlinearly Convergent Algorithm for A Class of Non-Convex Optimal Control Problems}, ArXiv 2112.
\newblock \href {http://arxiv.org/abs/2112.00108} {\path{arXiv:2112.00108}}.

\bibitem{Pavanello2024LT}
Z.~Pavanello, L.~Pirovano, R.~Armellin, A.~{De Vittori}, P.~{Di Lizia}, {Collision Avoidance Maneuver Optimization During Low-Thrust Propelled Trajectories}, Astrodynamics\href {http://dx.doi.org/https://doi.org/10.1007/s42064-024-0227-3} {\path{doi:https://doi.org/10.1007/s42064-024-0227-3}}.

\bibitem{DeVittori2022}
A.~{De Vittori}, M.~F. Palermo, P.~{Di Lizia}, R.~Armellin, {Low-Thrust Collision Avoidance Maneuver Optimization}, Journal of Guidance, Control, and Dynamics 45~(10) (2022) 1815--1829.
\newblock \href {http://dx.doi.org/10.2514/1.G006630} {\path{doi:10.2514/1.G006630}}.

\bibitem{Mueller2009}
J.~Mueller, {Onboard Planning of Collision Avoidance Maneuvers Using Robust Optimization}, in: AIAA Infotech@Aerospace Conference, American Institute of Aeronautics and Astronautics (AIAA), Reston, Virigina, 2009.
\newblock \href {http://dx.doi.org/10.2514/6.2009-2051} {\path{doi:10.2514/6.2009-2051}}.

\bibitem{Mueller2013}
J.~B. Mueller, P.~R. Griesemer, S.~J. Thomas, {Avoidance Maneuver Planning Incorporating Station-Keeping Constraints and Automatic Relaxation}, Journal of Aerospace Information Systems 10~(6) (2013) 306--322.
\newblock \href {http://dx.doi.org/10.2514/1.54971} {\path{doi:10.2514/1.54971}}.

\bibitem{DeVittori2025}
A.~{De Vittori}, Z.~Pavanello, P.~{Di Lizia}, J.~McMahon, R.~Armellin, {Combined Long-Term Collision Avoidance and Stochastic Station-Keeping in Geostationary Earth Orbit}, Journal of Guidance, Control, and Dynamics (2025) 1--15\href {http://dx.doi.org/10.2514/1.G008629} {\path{doi:10.2514/1.G008629}}.

\bibitem{Pavanello2023Long}
Z.~Pavanello, L.~Pirovano, R.~Armellin, {Long-Term Encounters Collision Avoidance Maneuver Optimization with a Convex Formulation}, in: 33rd AAS/AIAA Space Flight Mechanics Meeting, American Institute of Aeronautics and Astronautics (AIAA), Austin, TX, 2023, pp. 1--20.

\bibitem{Kim2012}
E.~H. Kim, H.~D. Kim, H.~J. Kim, {A Study on the collision avoidance maneuver optimization with multiple space debris}, Journal of Astronomy and Space Sciences 29~(1) (2012) 11--21.
\newblock \href {http://dx.doi.org/10.5140/JASS.2012.29.1.011} {\path{doi:10.5140/JASS.2012.29.1.011}}.

\bibitem{Morselli2014Collision}
A.~Morselli, R.~Armellin, P.~{Di Lizia}, F.~{Bernelli Zazzera}, F.~Bernelli-Zazzera, {Collision avoidance maneuver design based on multi-objective optimization}, in: AIAA/AAS Astrodynamics Specialist Conference 2014, San Diego, California, 2014, pp. 1819--1838.

\bibitem{Seong2015}
J.-D. Seong, H.-D. Kim, {Collision avoidance maneuvers for multiple threatening objects using heuristic algorithms}, Proceedings of the Institution of Mechanical Engineers, Part G: Journal of Aerospace Engineering 229~(2) (2015) 256--268.
\newblock \href {http://dx.doi.org/10.1177/0954410014530678} {\path{doi:10.1177/0954410014530678}}.

\bibitem{Seong2016}
J.~D. Seong, H.~D. Kim, {Multiobjective optimization for collision avoidance maneuver using a genetic algorithm}, Proceedings of the Institution of Mechanical Engineers, Part G: Journal of Aerospace Engineering 230~(8) (2016) 1438--1447.
\newblock \href {http://dx.doi.org/10.1177/0954410015611699} {\path{doi:10.1177/0954410015611699}}.

\bibitem{Masson2023}
M.~Masson, D.~Arzelier, M.~Joldes, B.~Revelin, \href{hal-03847541}{Multi-maneuvers algorithms for multi-risk collision avoidance via nonconvex quadratic optimization}, in: IFAC World Congress 2023, Yokohama, Japan, 2023.
\newline\urlprefix\url{hal-03847541}

\bibitem{Bourriez2023}
N.~Bourriez, A.~Loizeau, A.~F. Abdin, {Spacecraft autonomous decision-planning for collision avoidance: A Reinforcement Learning approach}, in: 74th International Astronautical Congress: IAC Proceedings, International Astronautical Federation (IAF), Baku, Azerbaijan, 2023, pp. 2--6.

\bibitem{Mu2024}
C.~Mu, S.~Liu, M.~Lu, Z.~Liu, L.~Cui, K.~Wang, {Autonomous spacecraft collision avoidance with a variable number of space debris based on safe reinforcement learning}, Aerospace Science and Technology 149~(4) (2024) 109131.
\newblock \href {http://dx.doi.org/10.1016/j.ast.2024.109131} {\path{doi:10.1016/j.ast.2024.109131}}.

\bibitem{Alfano2005Collision}
S.~Alfano, {Collision Avoidance Maneuver Planning Tools}, in: 15th AAS/AIAA Astrodynamics Specialist Conference, American Institute of Aeronautics and Astronautics (AIAA), Lake Tahoe, California, 2005.

\bibitem{Frigm2020}
R.~C. Frigm, M.~D. Hejduk, L.~C.~. Johnson, D.~Plakalovic, {Total Probability of Collision as a Metric for Finite Conjunction Assessment and Collision Risk Management}, in: Advanced Maui Optical and Space Surveillance Technologies (AMOS) Conference, Maui, Hawaii, 2015.

\bibitem{Bombardelli2014}
C.~Bombardelli, {Analytical formulation of impulsive collision avoidance dynamics}, Celestial Mechanics and Dynamical Astronomy 118~(2) (2014) 99--114.
\newblock \href {http://dx.doi.org/10.1007/s10569-013-9526-3} {\path{doi:10.1007/s10569-013-9526-3}}.

\bibitem{Bombardelli2015}
C.~Bombardelli, J.~Hernando-Ayuso, {Optimal impulsive collision avoidance in low earth orbit}, Journal of Guidance, Control, and Dynamics 38~(2) (2015) 217--225.
\newblock \href {http://dx.doi.org/10.2514/1.G000742} {\path{doi:10.2514/1.G000742}}.

\bibitem{Gonzalo2021}
J.~L. Gonzalo, C.~Colombo, P.~{Di Lizia}, {Analytical framework for space debris collision avoidance maneuver design}, Journal of Guidance, Control, and Dynamics 44~(3) (2021) 469--487.
\newblock \href {http://dx.doi.org/10.2514/1.G005398} {\path{doi:10.2514/1.G005398}}.

\bibitem{Abay2017}
R.~Abay, {Collision avoidance dynamics for optimal impulsive collision avoidance maneuvers}, in: 8th International Conference on Recent Advances in Space Technologies (RAST), IEEE, 2017, pp. 263--271.
\newblock \href {http://dx.doi.org/10.1109/RAST.2017.8002935} {\path{doi:10.1109/RAST.2017.8002935}}.

\bibitem{DeVittori2024capitolo6}
A.~{De Vittori}, {Impulsive Collision Avoidance Manoeuvres}, in: Enhanced Collision Avoidance Strategies in the Near-Earth Environment, Politecnico di Milano, Milan, IT, 2024, Ch.~2.

\bibitem{RubioAntón2024}
J.~Ant{\'{o}}n, F.~Biondi, D.~Escobar, P.~Di~Lizia, Fair shared collision avoidance manoeuvre for active vs active conjunctions, in: 75th International Astronautical Congress (IAC), IAF, Milano, Italy, 2024.

\bibitem{Wang2023}
Y.~Wang, F.~Topputo, Indirect optimization of fuel-optimal many-revolution low-thrust transfers with eclipses, IEEE Transactions on Aerospace and Electronic Systems 59~(1) (2023) 39--51.
\newblock \href {http://dx.doi.org/10.1109/TAES.2022.3189330} {\path{doi:10.1109/TAES.2022.3189330}}.

\bibitem{Lee2014}
K.~Lee, C.~Park, S.~Y. Park, {Near-optimal guidance and control for spacecraft collision avoidance maneuvers}, in: AIAA/AAS Astrodynamics Specialist Conference 2014, American Institute of Aeronautics and Astronautics (AIAA), Keystone, Colorado, 2014, pp. 1--20.
\newblock \href {http://dx.doi.org/10.2514/6.2014-4114} {\path{doi:10.2514/6.2014-4114}}.

\bibitem{Reiter2018}
J.~A. Reiter, D.~B. Spencer, {Solutions to rapid collision-avoidance maneuvers constrained by mission performance requirements}, Journal of Spacecraft and Rockets 55~(4) (2018) 1039--1047.
\newblock \href {http://dx.doi.org/10.2514/1.A33898} {\path{doi:10.2514/1.A33898}}.

\bibitem{Gonzalo2019Semi}
J.~L. Gonzalo, C.~Colombo, P.~{Di Lizia}, {A semi-analytical approach to low-thrust collision avoidance manoeuvre design}, in: 70th International Astronautical Congress (IAC), Vol.~70, Washington D.C., 2019, pp. 21--25.

\bibitem{Hernando-Ayuso2020}
J.~Hernando-Ayuso, C.~Bombardelli, {Low-Thrust Collision Avoidance in Circular Orbits}, Journal of Guidance, Control, and Dynamics 44~(5) (2021) 983--995.
\newblock \href {http://dx.doi.org/10.2514/1.G005547} {\path{doi:10.2514/1.G005547}}.

\bibitem{DeVittori2022Geo}
A.~{De Vittori}, G.~Dani, P.~{Di Lizia}, R.~Armellin, {Numerically Efficient Low-Thrust Collision Avoidance Maneuver Design in GEO Regime with Equinoctial Orbital Elements}, in: 2nd ESA NEO and Debris Detection Conference, no. January, ESA/ESOC, Darmstadt, Germany, 2023.

\bibitem{DeVittori2023}
A.~{De Vittori}, G.~Dani, P.~{Di Lizia}, R.~Armellin, {Low-Thrust Collision Avoidance Design for Leo Missions With Return To Nominal Orbit}, in: 33rd AAS/AIAA Space Flight Mechanics Meeting, American Institute of Aeronautics and Astronautics (AIAA), Austin, TX, 2023, pp. 1--18.

\bibitem{Strobel2023}
R.~Str\"{o}bel, E.~Stoll, Analytical assessment of short-duration low-thrust collision-avoidance maneuvers, Journal of Guidance, Control, and Dynamics 47~(2) (2024) 358--365.
\newblock \href {http://dx.doi.org/10.2514/1.G007175} {\path{doi:10.2514/1.G007175}}.

\bibitem{gonzalo2021computationally}
J.~L. Gonzalo, C.~Colombo, P.~Di~Lizia, et~al., Computationally efficient approaches for low-thrust collision avoidance activities, in: 72nd International Astronautical Congress: IAC Proceedings, International Astronautical Federation (IAF), Dubai, UAE, 2021, pp. 1--10.

\bibitem{gonzalo2022single}
J.~L. Gonzalo, C.~Colombo, P.~Di~Lizia, et~al., Single-averaged models for low-thrust collision avoidance under uncertainties, in: 73rd International Astronautical Congress: IAC Proceedings, International Astronautical Federation (IAF), Paris, France, 2022, pp. 1--8.

\bibitem{gonzalo37efficient}
J.~L. Gonzalo, C.~Colombo, P.~Di~Lizia, A.~De~Vittori, M.~Maestrini, P.~G. Padreny, M.~T. Ribell, D.~E. Ant{\'o}n, Efficient models for low thrust collision avoidance in space, Materials Research Proceedings 37.
\newblock \href {http://dx.doi.org/10.21741/9781644902813-137} {\path{doi:10.21741/9781644902813-137}}.

\bibitem{DeVittori2023a}
A.~{De Vittori}, M.~Omodei, P.~{Di Lizia}, R.~Armellin, P.~Gago, M.~{Torras Ribell}, J.~Ant{\'{o}}n, D.~Escobar, {Numerically Efficient Low-Thrust Fuel-Optimal Collision Avoidance Maneuvers With Tangential Firings}, in: AAS/AIAA Astrodynamics Specialist Conference, American Institute of Aeronautics and Astronautics (AIAA), Austin, Texas, 2023.

\bibitem{colombo2023sensitivity}
C.~Colombo, A.~De~Vittori, M.~Omodei, J.~L. Gonzalo, M.~Maestrini, P.~Di~Lizia, P.~Gago~Padreny, M.~Torras~Ribell, {\'A}.~Gallego~Torrego, D.~Escobar~Ant{\'o}n, et~al., Sensitivity analysis of collision avoidance manoeuvre with low thrust propulsion, in: Aerospace Europe Conference 2023-Joint 10th EUCASS-9th CEAS Conference, Lausanne, Switzerland, 2023, pp. 1--17.

\bibitem{DellElce2024}
L.~Dell'Elce, F.~{De Veld}, J.-B. Pomet, {Characterization of the minimum warning time for low-thrust collision avoidance maneuvers}, in: 2024 AAS/AIAA Astrodynamics Specialist Conference, American Institute of Aeronautics and Astronautics (AIAA), Broomfield, Colorado, 2024, pp. 1--19.

\bibitem{DeVeld2024}
F.~{De Veld}, D.~Lamberto, J.-B. Pomet, {Minimum Warning Time Analysis for Low-Thrust Collision Avoidance Manoeuvres with Steering Laws}, in: 75th International Astronautical Congress (IAC), Milan, Italy, 2024, pp. 1--8.

\bibitem{Polli2025}
E.~M. Polli, J.~L. Gonzalo, C.~Colombo, Semi-analytical model for autonomous fuel-optimal low-thrust collision avoidance with eclipse constraints, in: 2025 AAS/AIAA Space Flight Mechanics Meeting, American Institute of Aeronautics and Astronautics (AIAA), Kaua'i, Hawaii, 2025.

\bibitem{Zollo2024}
A.~Zollo, Z.~Pavanello, R.~Armellin, J.~F. {San Juan D{\'{i}}az}, B.~Schlepp, R.~Kahle, {Validation of a Fuel-Efficient Collision Avoidance Manoeuvre Optimiser for the GRACE-FO Mission}, in: 75th International Astronautical Congress (IAC), no. October, International Astronautical Federation (IAF), Milan, Italy, 2024, pp. 1--14.

\bibitem{DeMaria2023}
L.~{De Maria}, A.~{De Vittori}, P.~{Di Lizia}, {Numerically Efficient Impulsive and Low-Thrust Collision Avoidance Manoeuvres in Cislunar L1-Near Rectilinear Halo Orbit}, in: 74th International Astronautical Congress: IAC Proceedings, International Astronautical Federation (IAF), Baku, Azerbaijan, 2023.

\bibitem{Crum2024}
S.~Crum, M.~Borowitz, B.~C. Gunter, F.~Humphrey, \href{https://arc.aiaa.org/doi/abs/10.2514/6.2024-0843}{Cislunar orbit collision probability analysis}, in: AIAA SCITECH 2024 Forum, American Institute of Aeronautics and Astronautics (AIAA), 2024.
\newblock \href {http://dx.doi.org/10.2514/6.2024-0843} {\path{doi:10.2514/6.2024-0843}}.
\newline\urlprefix\url{https://arc.aiaa.org/doi/abs/10.2514/6.2024-0843}

\bibitem{lechtenberg2024novel}
T.~Lechtenberg, C.~J. Franz, J.~W. Gangestad, Novel strategies for cislunar conjunction assessment and collision avoidance, in: AIAA Aviation Forum and Ascend 2024, American Institute of Aeronautics and Astronautics (AIAA), 2024, p. 4821.

\bibitem{Nugent2025}
L.~Nugent, K.~Howell, D.~Davis, Flexible maneuver planning for collision avoidance in l1/l2 libration point orbits, in: 2025 AAS/AIAA Space Flight Mechanics Meeting, American Institute of Aeronautics and Astronautics (AIAA), Kaua'i, Hawaii, 2025.

\bibitem{Mishne2017}
D.~Mishne, E.~Edlerman, {Collision-Avoidance Maneuver of Satellites Using Drag and Solar Radiation Pressure}, Journal of Guidance, Control, and Dynamics 40~(5) (2017) 1191--1205.
\newblock \href {http://dx.doi.org/10.2514/1.G002376} {\path{doi:10.2514/1.G002376}}.

\bibitem{Omar2020}
S.~R. Omar, R.~Bevilacqua, {Spacecraft collision avoidance using aerodynamic drag}, Journal of Guidance, Control, and Dynamics 43~(3) (2020) 567--573.
\newblock \href {http://dx.doi.org/10.2514/1.G004518} {\path{doi:10.2514/1.G004518}}.

\bibitem{Turco2023}
F.~Turco, C.~Traub, S.~Gai{\ss}er, J.~Burgdorf, S.~Klinkner, S.~Fasoulas, {An analysis tool for collision avoidance manoeuvres using aerodynamic drag}, Acta Astronautica 211~(2) (2023) 116--129.
\newblock \href {http://dx.doi.org/10.1016/j.actaastro.2023.05.038} {\path{doi:10.1016/j.actaastro.2023.05.038}}.

\bibitem{gaglio2024optimal}
E.~Gaglio, C.~Traub, F.~Turco, J.~O. {Murcia Piñeros}, R.~Bevilacqua, S.~Fasoulas, Optimal drag-based collision avoidance: Balancing miss distance and orbital decay, Acta Astronautica 228 (2025) 295--305.
\newblock \href {http://dx.doi.org/https://doi.org/10.1016/j.actaastro.2024.11.052} {\path{doi:https://doi.org/10.1016/j.actaastro.2024.11.052}}.

\bibitem{Kleinig2022}
T.~Kleinig, B.~Smith, C.~Capon, {Collision avoidance of satellites using ionospheric drag}, Acta Astronautica 198~(April) (2022) 45--55.
\newblock \href {http://dx.doi.org/10.1016/j.actaastro.2022.03.017} {\path{doi:10.1016/j.actaastro.2022.03.017}}.

\end{thebibliography}

\end{document}